\documentclass[10pt]{article}
\usepackage{fancyvrb}  
\usepackage{amsmath,amssymb,color,bbm,ifthen,enumerate}
\usepackage{graphicx}
\usepackage{subfigure}
\usepackage{comment}

\newtheorem{thm}{Theorem}[section]

\newtheorem{lem}[thm]{Lemma}

\newtheorem{cor}[thm]{Corollary}

\newenvironment{pf}{\paragraph{Proof.}}
{\nopagebreak\hfill\nopagebreak\rule{2mm}{2mm}\par\bigskip}

\newcommand{\bone}{\mbox{\boldmath$1$}}

\newcommand{\pt}{\frac{\partial}{\partial t}}

\newcommand{\prob}{\mathrm{Prob}}
\newcommand{\krnl}{g}
\newcommand{\cdf}{F}

\begin{document}

\title{Distribution functions of Poisson random integrals: \\ Analysis and computation}

\author{ Mark Veillette  and Murad S. Taqqu\thanks{ This research was partially supported by the NSF grants
DMS-0706786, and DGE-0221680.}
\\
Boston University 
}


\maketitle

\begin{abstract}

We want to compute the cumulative distribution function of a
one-dimensional Poisson stochastic integral $I(\krnl) = \displaystyle
\int_0^T \krnl(s) N(ds)$, where $N$ is a Poisson random measure with
control measure $n$ and $\krnl$ is a suitable kernel function.  We do
so by combining a Kolmogorov-Feller equation with a finite-difference
scheme.  We provide the rate of convergence of our numerical scheme
and illustrate our method on a number of examples.   The software used
to implement the procedure is available on demand and we demonstrate
its use in the paper.  
    
\end{abstract}

\renewcommand{\baselinestretch}{1.3}

\section{Introduction}

Let $T>0$, and let $N(\cdot)$ be a 
Poisson random measure defined on the interval $[0,T]$, with the Borel
$\sigma$-field, and control
measure $n(ds)$, which we assume to have a density $n(ds) = n(s)ds$.  For appropriate functions $\krnl(s), \ s \in [0,T]$, the Poisson stochastic integral
\begin{equation}\label{e:firstpoisint}
I(\krnl)  =  \int_0^T \krnl(s) N(ds),
\end{equation}
is a random variable defined as the limit in probability:
\begin{equation}
P-\lim_{||\Delta_n|| 
\rightarrow 0} \sum_{i=1}^n \krnl(s_i) N( [s_i,s_{i+1}) ), \quad \Delta_n
= \{ 0 < s_1 < s_2 < \dots <
s_n < T \}.
\end{equation}
For this limit to exist, we require that kernel $\krnl$ satisfies 
\begin{equation} \label{e:poiscond}
\int_0^T \min(\krnl(s),1) n(ds) < \infty.
\end{equation} 
The characteristic function of $I(\krnl)$ is expressible in terms
of the control measure and kernel, and is given by
\begin{equation}\label{e:poischarfunc}
\phi_\krnl(\theta) = \mathbb{E} e^{i \theta I(\krnl) } = \exp \left( \int_0^T (e^{i \theta \krnl(s)}
  - 1) n(ds) \right). 
\end{equation}
For more details about Poisson integrals and some discussion of
applications, see \cite{kingman:1993}, Chapter 3.  

It is important to note that in the case where $n$ is a finite
measure on $(0,T)$, the distribution of $I(g)$
is that of a compound Poisson distribution.  This can be seen in the case
where $g$ is strictly increasing and has inverse $g^{-1}$ by making a change of variable $u = g(s)$ in the integral in
(\ref{e:poischarfunc}):
\begin{equation*}
\exp \left( \int_0^T (e^{i \theta \krnl(s)}
  - 1) n(ds) \right) = \exp\left( \int_{g^{-1}(0)}^{g^{-1}(T)} (e^{i
    \theta u} - 1) \bar{n}(u) du \right)
\end{equation*}
where $\bar{n}(u) =  (g^{-1})'(u) n(g^{-1}(u))$, which is the
characteristic function of a compound Poisson distribution
(see Proposition 1.2.11 in \cite{Applebaum:2004}).  Thus, $I(g)$ in this
case has a compound Poisson distribution with rate $r = \int_{g^{-1}(0)}^{g^{-1}(T)} \bar{n}(u) du = \int_0^T n(s) ds$
and jump distribution with density $r^{-1} \bar{n}(u)$, $g^{-1}(0) < u
< g^{-1}(T)$.  A similar argument
shows this for non-monotone $g$, except one must construct $\bar{n}$ carefully by
breaking the integral (\ref{e:firstpoisint}) into pieces on which $g$
is either flat or strictly monotone.   Thus, by considering stochastic
integrals of this type, we are covering a large class of interesting
distributions.  

The distribution of integrals such as (\ref{e:firstpoisint}) would be
easy to obtain if the measure $N$ was Gaussian.  In this case $I(\krnl)$ is also Gaussian: $I(\krnl) \sim N(0,\int_0^T \krnl^2(s) n(s) ds )$.  Since $N$ is Poisson, however, the integral $I(\krnl)$ is not in general Poisson and its distribution is not easy to determine.  Our goal here is to
 study the cumulative distribution function (CDF) of $I(\krnl)$, which we denote by $\cdf(x) = \mathrm{Prob}(I(\krnl) \leq x)$ and develop a convenient numerical scheme to evaluate it.  In particular, we focus on the following:
\begin{itemize}
\item[1.] A Kolmogorov-Feller type evolution equation associated to $\cdf$.
\item[2.] Smoothness properties of $\cdf$.  More specifically, we show that under certain assumptions on $\krnl$ and $n$, $\cdf$ lives in the H\"{o}lder space $C^{0,\gamma}((\krnl(0),\krnl(T))$, for some $0<\gamma \leq 1$ depending on $\krnl$.
\item[3.] A numerical method for computing $\cdf$.  
\end{itemize}

At first glance, one might attempt to compute the CDF $\cdf$ using the characteristic function (\ref{e:poischarfunc}) and the integration-based inversion methods of Abate and Whitt \cite{abate:1992}.   For continuous $\cdf$, their approach boils down to computing the following integral numerically:
\begin{equation}\label{e:abatesmet}
\cdf(x) = \frac{ 2 }{ \pi } \int_0^\infty \mathrm{Re}(\phi_\krnl)(u) \frac{\sin x u}{u} du.
\end{equation}
  However, this method is not always efficient in this setting for
  multiple reasons. First, evaluating the characteristic function
  $\phi_\krnl$ is not always easy, since $\int_0^T (e^{i \omega
    \krnl(s)} - 1)n(ds)$ is unlikely to have a closed form, and
  numerically computing this integral might be difficult.  This adds
  another source of error on top of the ``truncation'' and
  ``discretization'' errors associated with the integration of
  (\ref{e:abatesmet}).  Moreover, these errors are difficult to bound
  exactly.   Secondly, we will see in the following that $\cdf$ is not
  differentiable in general, which translates to a slow rate of decay
  of $\phi_\krnl$.  This makes numerical integration difficult.  

We propose an alternative method for computing $\cdf$ which does not involve (directly) the characteristic function (\ref{e:poischarfunc}).  Our method has several advantages over brute-force integration of (\ref{e:abatesmet}).  First, our method only requires knowing $\krnl$ and the density of the control measure $n$, and thus is much faster than integration methods which must evaluate the characteristic function (\ref{e:poischarfunc}).  Second, we provide here error bounds given general assumptions on $\krnl$ and $n$.  Third, our method generates the CDF $\cdf$ on an entire interval instead of a single point.  Fourth, we obtain the more general cumulative distribution functions $\cdf(x,t)$, $t \geq 0$, of the following {\it stochastic process}:
\begin{equation}\label{e:xzeroplusint}
X(t) = X_0 + \int_0^t \krnl(s) N(ds), \quad 0 \leq t \leq T,
\end{equation} 
where $X_0$ is an independent random variable with given CDF
$\cdf_0(x)$.  In the sequel, we illustrate some of these advantages on
an example.

Our method involves solving an evolution equation satisfied by the CDF of $X(t)$ with initial condition $\cdf_0$. A nice way to view such an evolution equation is to
consider the characteristic function $\phi_\krnl(\omega,t)$ of (\ref{e:xzeroplusint})  (we will from hereon drop the ``$\krnl$'' subscript from $\phi$).  Using (\ref{e:poischarfunc}), we have at time $0 \leq t
\leq T$,
\begin{equation}
\phi(\theta,t) = \phi_0(\theta)  \exp \left( \int_0^t (e^{i \theta \krnl(s)}
  - 1) n(ds) \right),
\end{equation}
where $\phi_0$ is the characteristic function of $X_0$.     
Differentiating with respect to $t$ gives
\begin{equation}\label{e:poisftODE}
\pt \phi(\theta,t) = \left( (e^{i \theta \krnl(t)}
  - 1) n(t)  \right) \phi(\theta,t).
\end{equation}
Thus, the characteristic function $\phi(\omega,t)$ can be viewed as the solution
at time $t$ of the simple ordinary differential equation
(\ref{e:poisftODE}) with initial condition $\phi(\theta,0) =
\phi_0(\theta)$.  With this ODE perspective, the naive approach of brute-force integration can be thought of as a two step process: (i) evolving $\phi_0$
under the trivial dynamics of (\ref{e:poisftODE}) to obtain $\phi(\theta,t)$, and then (ii) computing (\ref{e:abatesmet}) to get $\cdf$. 

An alternative approach of obtaining an evolution equation for $\cdf$ is to notice that the process $X$ in (\ref{e:xzeroplusint}) is a continuous time Markov process, hence it satisfies the {\it Kolmogorov-Feller forward equation}  (\cite{kannan:1979}, Section 5.1).  In the infinitesimal time interval $[t,t+dt)$, $X$ may take a jump of size $\krnl(t)$ with probability $n(t) dt + o(dt)$,  thus, the Kolmogorov-Feller equation for the single-time density $p(x,t)$ of the process $X$ takes the form 
\begin{equation}\label{e:genKol}
\frac{\partial}{\partial t} p(x,t)  = - n(t) (p(x,t) - p(x - \krnl(t) ,t)), \quad t \geq 0, \ x \in \mathbb{R}.
\end{equation}
Integrating this equation suggests that a similar equation should hold for $\cdf$, however some technicalities arise because $p$ can contain ``atoms''.   We show directly that a similar equation holds for the CDF, with the exception of points where the CDF is discontinuous.  We then solve this evolution equation numerically using a finite-difference scheme.  

This paper is outlined as follows.  In Section \ref{s:mainthm}, we
state a theorem which gives an evolution equation of the function $\cdf$.  We then study the smoothness of the function $\cdf(\cdot,t)$ in Section \ref{s:holder}.  We then describe in Section \ref{s:methods} a numerical method for calculating
$\cdf$ using a finite-difference scheme, and provide a rate of convergence for our method.  In Section \ref{s:examples}, we apply our methods to a
collection of examples and in Section \ref{s:approx}, we
establish lemmas related to the quality of the approximation.
Finally, in Section \ref{s:software} we provide a guide to the
software, written in MATLAB, which allows the user to obtain the cumulative distribution function and density
function of the integral $\int_0^T
\krnl(s) N(ds)$ for
general $\krnl$.

\section{Kolmogorov-Feller equation for $\cdf$}\label{s:mainthm}

In this section we show that the CDFs $\cdf(x,t)$, $ t \geq 0$ of $X(t)$ satisfy a Kolmogorov-Feller-type equation.  As mentioned in the introduction, the form of this equation is easy to guess based on the Kolmogorov-Feller equation  (\ref{e:genKol}) satisfied by the density of the process $X(t)$.  

Since we will be dealing with distribution functions, we must be
mindful of discontinuities which may arise. Given a CDF $\cdf$, we will let $C(\cdf)$ denote the points of continuity
of $\cdf$. Also, we define the \textit{right}-time derivative as
\begin{equation}
\frac{\partial}{\partial t^+} \cdf(x,t) = \lim_{h \rightarrow 0+}
\frac{\cdf(x,t+h) - \cdf(x,t)}{h}.
\end{equation}   

The following theorem specifies the equation for $\cdf$ which will be used the sequel. 

\begin{thm}\label{t:mainthm}
Let $\krnl(s), \ s \geq 0$ be a right-continuous function satisfying (\ref{e:poiscond}), and consider the process
\begin{equation}
X(t) = X_0 + \int_0^t \krnl(s) N(ds), \quad 0 \leq t \leq T,
\end{equation} 
where $N$ is a Poisson random measure with control measure $n(ds) =
n(s)ds$ and $X_0$ is a given random variable.  Let $\cdf(x,t)$ be the CDF
of $X$ at time $t$.  Then, for all $(x,t)$ such that $x-\krnl(t) \in C(\cdf(\cdot,t))$, $\cdf$ satisfies the
following equation:
\begin{equation}\label{e:mainde}
\frac{\partial}{\partial t^+} \cdf(x,t)  = -n(t)(\cdf(x,t) - \cdf(x-\krnl(t),t)).
\end{equation}
\end{thm}

\begin{pf}

Let $x,t \geq 0$ such that $x - \krnl(t) \in C(\cdf)$ and let $h>0$.  Consider the difference 
\begin{equation}
\cdf(x,t+h) - \cdf(x,t) = \mathbb{E}[ \bone_{(-\infty, x]}(X(t+h)) -  \bone_{(-\infty,x]}(X(t))].
\end{equation}
By conditioning
on the number $N$ of Poisson arrivals in the interval $[t,t+h)$, the expectation on the right-hand side becomes
\begin{equation*}
\mathbb{E}[ \bone_{(-\infty, x]}(X(t+h)) -  \bone_{(-\infty,x]}(X(t))]
= \sum_{n=0}^\infty  \mathbb{E}[ \bone_{(-\infty, x]}(X(t+h)) -
\bone_{(-\infty,x]}(X(t)) | N=n]\ \prob[N=n] 
\end{equation*}
\begin{equation}\label{e:cond}
=  \mathbb{E}[ \bone_{(-\infty, x)}(X(t+h)) -
\bone_{(-\infty,x]}(X(t)) | N=1]  \left( \int_t^{t+h} n(\tau)d\tau \right)
e^{-\int_t^{t+h} n(\tau)d\tau} + O(h^2).
\end{equation}
Here we have used that facts that if $N=0$, $X(t)  = X(t+h)$, and the
probability of seeing two or more Poisson arrivals in the interval
$[t,t+h)$ is $O(h^2)$ since we've assumed $n$ is continuous.  In the
event of one arrival in $[t,t+h)$, namely $N=1$, we have $X(t+h) = X(t) + \krnl(S)$,
where $S \in [t,t+h)$ is a random time with density
\begin{equation}   
S \sim f_S(s) =  \begin{cases} \displaystyle \frac{n(s)}{\int_t^{t+h} n(\tau) d\tau}, \quad s
  \in [t,t+h) \\ 0, \qquad \qquad \qquad \mathrm{else} \end{cases}.
\end{equation}
Conditioning on $S$,  (\ref{e:cond}) now becomes
\begin{equation*}
\left( \int_t^{t+h} \mathbb{E}[ \bone_{(-\infty, x]}(X(t)+\krnl(s)) -
\bone_{(-\infty,x]}(X(t)) \ | \ S = s]  f_S(s) ds \right) \left( \int_t^{t+h} n(\tau)d\tau \right)
e^{-\int_t^{t+h} n(\tau)d\tau} + O(h^2)
\end{equation*}
\begin{equation*}
= \left( \int_t^{t+h} \mathbb{E}[ \bone_{(-\infty, x]}(X(t)+\krnl(s)) -
\bone_{(-\infty,x]}(X(t))] n(s) ds \right)
e^{-\int_t^{t+h} n(\tau)d\tau} + O(h^2)
\end{equation*}
\begin{equation*}
 =\left( \int_t^{t+h} \mathbb{E}[ \bone_{(-\infty, x-\krnl(s)]}(X(t)) -
\bone_{(-\infty,x]}(X(t))] n(s) ds \right)
e^{-\int_t^{t+h} n(\tau)d\tau} + O(h^2)
\end{equation*}
\begin{equation*}
= \left( \int_t^{t+h} (\cdf(x-\krnl(s),t) - \cdf(x,t) )n(s) ds \right)
e^{-\int_t^{t+h} n(\tau)d\tau} + O(h^2).
\end{equation*}
Now, since $\cdf(\cdot,t)$ is continuous at $x-\krnl(t)$, $\krnl$ is right continuous,
and $n$ is continuous,  let $h \rightarrow 0^{+}$ and get
\begin{eqnarray}
\lim_{h\rightarrow 0^+} \frac{\cdf(x,t+h) - \cdf(x,t)}{h} &=& \lim_{h
  \rightarrow 0^+}  \left( h^{-1} \int_t^{t+h} (\cdf(x-\krnl(s),t) - \cdf(x,t) )n(s) ds \right)
e^{-\int_t^{t+h} n(\tau)d\tau}   \label{e:takelim} \\
&=& n(t) (\cdf(x - \krnl(t),t) - \cdf(x,t)),
\end{eqnarray}
  which completes the proof.

\end{pf}

{\bf Remarks.}
\begin{itemize}

\item[1.] If $\cdf$ and $\krnl$ are continuous, then the right and left time derivatives
of $\cdf$ coincide and the left hand side of
(\ref{e:mainde}) can be replaced with $\pt \cdf(x,t)$.

\item[2.] If $ x - g(t) \notin C(\cdf)$ and $\krnl \geq 0$ is continuous and monotone increasing, then $ x- \krnl(s) \nearrow x - \krnl(t)$ as $s \searrow t$ and the limit in (\ref{e:takelim}) still exists, except we have instead
\begin{equation}\label{e:leftversion}
\frac{\partial \cdf}{\partial t^{+}} = n(t) ( \cdf^{-}(x-g(t),t) - \cdf(x,t)),
\end{equation}
where $\cdf^{-}$ is the left-continuous version of $\cdf$, i.e.~$\cdf^{-}(x,t) = \lim_{y \rightarrow x^- }\cdf(y,t)$.  This fact will be used in the sequel.

\item[3.] The continuity set $C(\cdf)$ can be identified from $\krnl$ and $n$ as follows.  We can write 
\begin{equation}
I(\krnl) = X_0 +  \int_0^t \krnl(s) N(ds) = X_0 + \sum_{i=1}^N \krnl(S_i),
\end{equation}
where $N$ is a Poisson random variable with mean $\int_0^t n(s) ds$.  Then, the atoms of $I(\krnl)$ correspond to the atoms of the random sums $ X_0 + \sum_{i=1}^n \krnl(S_i)$, with $n=0,1,2,\dots$ and $S_i$ are the ``Poisson events'', which are i.i.d.~with density 
\begin{equation}
\frac{n(s)}{\int_0^t n(\tau) d \tau}, \quad 0 < s < t.
\end{equation}
The continuity set $C(\cdf)$ is the complement of this set of atoms.  

For example, if $X_0 = 0$, $\krnl \equiv \eta$ is a constant function, and $n(s) \equiv 1$, then the atoms of the random sums $\sum_{i=1}^N g(S_i)$, which is the set $\{n\eta,n=0,1,2,\dots \}$.

\item[4.] To arrive at Theorem \ref{t:mainthm}, we used properties of the Poisson integral to derive the evolution equation (\ref{e:mainde}).  This reasoning can also be turned on its head -- this theorem can give a probabilistic interpretation to initial value problems for a class of differential-difference equations of the form 
\begin{equation}\label{e:PoisIVP}
\begin{cases}  \displaystyle\frac{\partial u}{\partial t^{+}}  = - n(s) ( u(x,t) - u(x - \krnl(t) , t ) ) , \quad x \in \mathbb{R}, t \geq 0 \\
u(x,0) = u_0(x) \end{cases},
\end{equation}
where $n$, $\krnl$ and $u_0$ are given functions and $u$ is an unknown function of $x$ and $t$.  By computing the CDF of the Poisson stochastic integral (\ref{e:firstpoisint}), we are essentially computing the {\it fundamental solution} of this differential-difference equation, since the solution to (\ref{e:PoisIVP}) can be written as the convolution of $u_0$ with the CDF $\cdf(x,t)$ of the integral (\ref{e:firstpoisint}) up to time $t$:
\begin{equation}\label{e:ddtruesol}
 u(x,t) = \int_{-\infty}^\infty u_0(x - u) d\cdf(u,t).
\end{equation}

To illustrate this point, consider the following example:  Let $b(t), \ t \geq 0$ be a given  smooth function with derivative $b'(t)$, and recall the classical transport PDE:
\begin{equation}\label{e:transport}
\begin{cases}  \displaystyle \frac{\partial u}{\partial t} = -b'(t) \frac{\partial u}{\partial x}, \quad x \in \mathbb{R}, t \geq 0 \\
u(x,0) = u_0(x)  \end{cases},
\end{equation}
where, for simplicity, the initial condition $u_0$ is some smooth and bounded function.  The solution to this equation is given by 
\begin{equation}\label{e:soldrift}
u(x,t) = u_0( x - b(t) ), \quad t \geq 0.
\end{equation}  
Now, let $h > 0$, and say we want to interpret the solution to the following variation of (\ref{e:transport}) where we replace $\partial u/\partial x$ with a finite-difference approximation:
\begin{equation}\label{e:transport1}
\begin{cases}  \displaystyle \frac{\partial u}{\partial t} = -b'(t) \frac{ u(x,t) - u(x - h,t)  }{ h }, \quad x \in \mathbb{R}, t \geq 0 \\
u(x,0) = u_0(x)  \end{cases}.
\end{equation}    
Notice that (\ref{e:transport1}) is of the form (\ref{e:mainde}) with $\krnl(t) \equiv h$, and $n(t) = b'(t)/h$.  This corresponds to the trivial Poisson random integral
\begin{equation}
\int_0^t h N(ds), \quad n(s) = \frac{b'(s)}{h},
\end{equation}
which has the distribution of the random variable $h N([0,t])$, where $N([0,t])$ has a Poisson distribution with mean $ b(t) h^{-1} $.    

Observe that in this case, (\ref{e:ddtruesol}) can be expressed as $u(x,t) = \mathbb{E}u_0(x - hN([0,t]))$ which is the convolution of $u_0$ with the CDF of $h N([0,1])$. Compare this to the solution of the classical transport equation (\ref{e:soldrift}).   

\end{itemize}

\section{Smoothness of $\cdf(\cdot,t)$}\label{s:holder}

Our goal is to develop a numerical scheme for approximating the CDF $\cdf$ which is based on the differential-difference equation (\ref{e:mainde}).  In order to estimate the error in our approximation to the true CDF, we must have some smoothness properties for $\cdf$, and thus make some assumptions on $\krnl$.  We will assume that $\krnl$ has a continuous inverse.  While this seems like a restrictive assumption, we discuss ways of generalizing our method to non-invertible $\krnl$ in section (\ref{s:extend}). 

\begin{figure}[ht]
\centering
\includegraphics[width=1.0\textwidth]{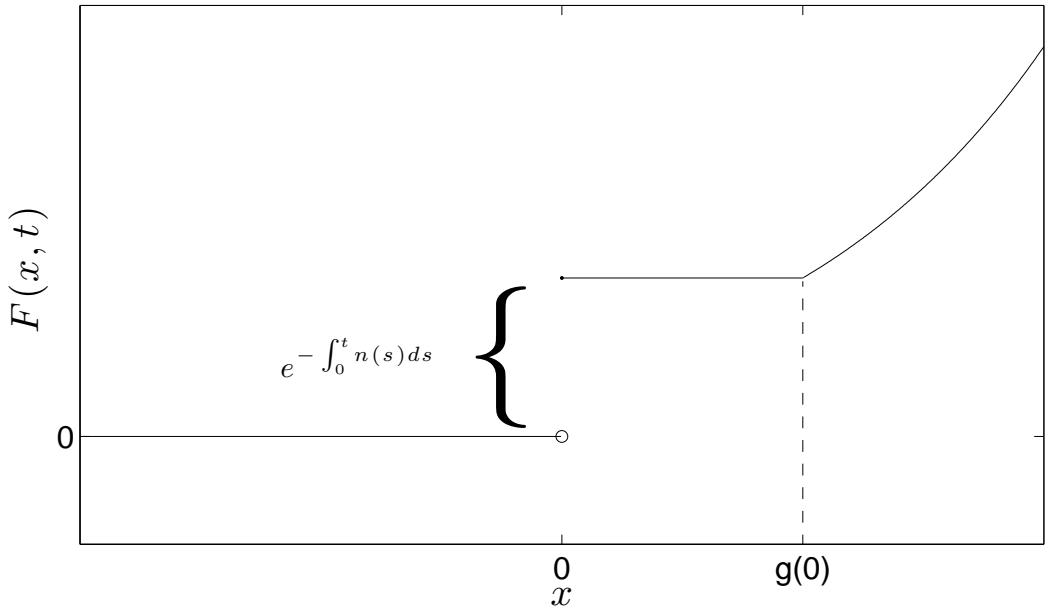}
\caption{ This shows the behavior of $\cdf(x,t)$ near $x = 0$ with $t$ fixed.  We see a jump of size $\exp(-\int_0^t n(s) ds )$ at $x=0$ which equals the probability of zero Poisson events in time $(0,t)$. Also, $\cdf(x,t)$ is constant in the interval $[0,\krnl(0))$ since with one or more Poisson event in $(0,t)$, $\krnl(0)$ is the minimum value the integral can take.  }
\label{f:pcartoon}
\end{figure}

In the following, we assume that $\krnl$ is a non-negative, strictly monotone continuous function on $[0,t]$ with inverse $\krnl^{-1}$.   We will also suppose that $n(s)$ is a  bounded function.   Since $\krnl \geq 0$, the integral $\int_0^t \krnl(s)N(ds)$ is zero if there are no Poisson events in $[0,t)$, i.e.~if $N[0,t) =0$ which happens with probability $e^{-\int_0^t n(s) ds}$.  Also,  $\int_0^t \krnl(s)N(ds) \geq \krnl(0)$ if $N[0,t) \geq 1$.  Therefore, $\cdf(x,t) = 0$ for $x < 0$, $\cdf$ has a jump of size $e^{-\int_0^t n(s)ds }$ at $x = 0$, and $\cdf \equiv e^{-\int_0^t n(s)ds }$ on the interval $[0,\krnl(0))$ (see Figure \ref{f:pcartoon}).  Having identified the discontinuity at $x=0$ in the CDF of $\cdf$, we now focus on the following question: How smooth is $\cdf$ for $x>0$?  This will depend on $\krnl$. 

Let $U \subset \mathbb{R}$.  For functions $u$ defined on $ U$ and for $0 < \gamma \leq 1$, recall the definition of the H\"{o}lder seminorm $[ \cdot ]_{C^{0,\gamma}(U) }$ and the H\"{o}lder space $C^{0,\gamma}( U )$ defined as
\begin{eqnarray}
[u]_{C^{0,\gamma} (U) } &=& \sup_{ \underset{x\neq y } { x,y \in U  }}  \frac{ |u(x) - u(y)| }{ |x - y|^\gamma } \\
C^{0,\gamma}(U) &=& \{ u\ : \ [u]_{C^{0,\gamma}(U) } < \infty \}.
\end{eqnarray}
The exponent $\gamma$ is called the H\"{o}lder exponent of $u$. For more on these spaces, see \cite{evans:1998} section 5.1.

The next theorem shows that $\cdf(x,t)$ for $x>0$ and $t>0$ lies in the same H\"{o}lder space as $\krnl^{-1}$.

\begin{thm}\label{t:holder}
Let $\krnl$ be a  non-negative, strictly increasing continuous function on $[0,T]$ with inverse $\krnl^{-1}$. For $0<t < T$ fixed, let $\cdf(x,t) = \mathrm{Prob}[X(t) \leq x]$, where $X(t)$ is defined by
\begin{equation}
X(t) =  \int_0^t \krnl(s) N(ds),
\end{equation}
where $N$ is a Poisson random measure with control measure $n(s)ds$ with $n$ a bounded function.  Let 
\begin{equation}
n^{\star}_t = \sup_{0 \leq s \leq t} n(s).
\end{equation}
  Then, if $\krnl^{-1} \in C^{0,\gamma}( ( \krnl(0),\krnl(t) ) ) $ with $0 < \gamma \leq 1$, then $\cdf(\cdot,t) \in C^{0,\gamma}( (0,\infty) )$ and
\begin{equation}
[ \cdf(\cdot,t) ]_{C^{0,\gamma} (0,\infty)} \leq C_{n,t} \ [\krnl^{-1}]_{C^{0,\gamma} (\krnl(0),\krnl(t))} ,
\end{equation}
where the constant $C_{n,t}$ is given by 
\begin{equation}
C_{n,t} =  \frac{ (1 - e^{\int_0^t n(s)ds})}{\int_0^t n(s) ds} n^{\star}_t
\end{equation}

\end{thm}

\begin{pf}


Fix $t>0$ and let $S,S_1,S_2\dots$ be i.i.d.~random variables with density 
\begin{equation}\label{e:defofnbar}
n_t(s) = \frac{n(s)}{\int_0^t n(\tau) d \tau}, \quad 0 < s < t.
\end{equation}
Let $\bar{n}_t(x) = \int_0^x n_t(s) ds $ be the CDF of $S$ and let 
\begin{eqnarray}
G(x) = \mathrm{Prob}[\krnl(S) \leq x] &=& \begin{cases}  0 \quad x < \krnl(0)  \\
\mathrm{Prob}[S < \krnl^{-1}(x)] \quad \krnl(0) \leq x \leq \krnl(t) \\
1 \quad  x > \krnl(t) \end{cases} \nonumber  \\ 
&=& \begin{cases} 0 \quad x < \krnl(0)  \\
\bar{n}_t(\krnl^{-1}(x)) \quad \krnl(0) \leq x \leq \krnl(t) \\
1 \quad  x > \krnl(t) \end{cases}.  \label{e:defofGnbar}
\end{eqnarray}
By conditioning on the number of arrivals $N$ in the interval $[0,t]$, we have for $ x,y > 0 $,
\begin{equation}
|\cdf(x,t) - \cdf(y,t)| = \left| \sum_{n=1}^\infty \left( \mathrm{Prob}[\sum_{i=1}^n \krnl(S_i) \leq x] - \mathrm{Prob}[\sum_{i=1}^n \krnl(S_i) \leq y] \right) P_n\right|
\end{equation}
where $P_n = \mathrm{Prob}[N=n]$.  Let $G^{\star n}$ denote the n-fold convolution of $G$ with itself. Now we have
\begin{eqnarray}
|\cdf(x,t) - \cdf(y,t)| &\leq& | G(x) - G(y) | P_1 + \sum_{n=2}^\infty | G^{\star n}(x) - G^{\star n}(y) | P_n \nonumber \\
&=&  | G(x) - G(y) | P_1 + \sum_{n=2}^\infty \left| \int_0^x G(x-u) G^{\star (n-1)}(du) - \int_0^y G(y-u) G^{\star (n-1)}(du) \right| P_n \nonumber \\
&\leq&  | G(x) - G(y) | P_1 + \sum_{n=2}^\infty \int_0^{\max(x,y)} |G(x-u) -G(y-u)| G^{\star (n-1)}(du) P_n 
\end{eqnarray}
Thus, for $0 < \gamma \leq 1$,
\begin{eqnarray}
\frac{|\cdf(x,t) - \cdf(y,t)| }{ |x - y|^\gamma } &\leq& \frac{| G(x) - G(y) | }{ |x - y|^\gamma } P_1 + \sum_{n=2}^\infty \int_0^{\max(x,y)} \frac{| G(x-u) - G(y-u) | }{ |x - y|^\gamma } G^{\star (n-1)}(du) P_n  \\
&\leq&  [G]_{ C^{0,\gamma} ((0,\infty)) } \sum_{n=1}^\infty P_n \\
&=&  ( 1 - e^{-\int_0^t n(s) ds })  [G]_{ C^{0,\gamma} ((0,\infty)) } \label{e:bound1}
\end{eqnarray}

What remains is to show that $ [G]_{ C^{0,\gamma} ((0,\infty)) } $ is bounded by a constant multiple of $[\krnl^{-1}]_{ C^{0,\gamma}(\krnl(0),\krnl(t)) }$.  Observe that $\bar{n}'_t(x) = n_t(x)$ where $n_t$ is defined in (\ref{e:defofnbar}) and thus
 $|\bar{n}_t'(x)| \leq n^{\star}_t \cdot  (\int_0^t n(s) ds)^{-1}$, where $n_t^{\star} = \sup_{0 \leq s \leq t} n(s)$. In view of (\ref{e:defofGnbar}),  the mean value theorem implies that for $\krnl(0) \leq x,y \leq \krnl(t)$,
\begin{equation}
\frac{|G(x) - G(y)|}{|x-y|^\gamma} =\frac{ |\bar{n}_t(\krnl^{-1}(x)) -  \bar{n}_t(\krnl^{-1}(y))| }{|x-y|^\gamma}  \leq \frac{n^{\star}_t}{\int_0^t n(s) ds} \frac{ |\krnl^{-1}(x) -  \krnl^{-1}(y)| }{|x-y|^\gamma} \leq \frac{n^{\star}_t}{\int_0^t n(s) ds}  [\krnl^{-1}]_{ C^{0,\gamma}(\krnl(0),\krnl(t)) } .
\end{equation}
If $0 < x < \krnl(t) \leq y$, then $G(y) = 1 = \bar{n}_t(t) = \bar{n}_t( \krnl^{-1}(\krnl(t)))$ and
\begin{equation}
\frac{|G(x) - G(y)|}{|x-y|^\gamma} \leq \frac{|\bar{n}_t(\krnl^{-1}(x)) - \bar{n}_t(\krnl^{-1}(\krnl(t)))   |}{|x-\krnl(t)|^\gamma} \leq   \frac{n^{\star}_t}{\int_0^t n(s) ds}   [\krnl^{-1}]_{ C^{0,\gamma}(\krnl(0),\krnl(t)) }. 
\end{equation}
A similar argument shows this also holds for $x < \krnl(0) < y$.  Finally, if $x,y \geq \krnl(t)$ or $x,y \leq \krnl(0)$, then 
\begin{equation*}
G(x) - G(y)  = 0 \leq  \frac{n^{\star}_t}{\int_0^t n(s) ds}  [\krnl^{-1}]_{ C^{0,\gamma}(\krnl(0),\krnl(t)) }.  
\end{equation*}
Thus, in all cases,
\begin{equation}\label{e:bound2}
\frac{|G(x) - G(y)|}{|x-y|^\gamma}  \leq   \frac{n^{\star}_t}{\int_0^t n(s) ds}   [\krnl^{-1}]_{ C^{0,\gamma}(\krnl(0),\krnl(t)) }.
\end{equation}
Hence,  (\ref{e:bound1}) and (\ref{e:bound2}) imply
\begin{equation} 
[\cdf(\cdot,t)]_{C^{0,\gamma}((0,\infty))}  \frac{ n^{\star}_t ( 1 - e^{-\int_0^t n(s) ds })}{\int_0^t n(s) ds}  [\krnl^{-1}]_{ C^{0,\gamma}(\krnl(0),\krnl(t)) } = C_{n,t} \ [\krnl^{-1}]_{ C^{0,\gamma}(\krnl(0),\krnl(t)) }  .
\end{equation}
which finishes the proof.

\end{pf}

In the preceding Theorem, we assumed $g$ is monotone.  What happens if $\krnl$ is piecewise monotone, i.e.~$\krnl(t)$ is monotone on the intervals $[t_i,t_{i+1})$, for $i=0,1,\dots,M$?  Then, the integral of $\krnl$ over $[t_0,t_M)$ can be ``constructed'' inductively by considering the sums
\begin{equation}   \label{e:poissummands}
I_i = I_{i-1} + \int_{t_{i-1}}^{t_i} \krnl(s) N(ds), \quad i=1,2,\dots,M
\end{equation}
with $I_0 = 0$.  Since $N$ is a Poisson measure, the summands in (\ref{e:poissummands})  are independent random variables.  Hence, the full CDF of $\int_{t_0}^{t_M}\krnl(s)N(ds)$ can be seen as the successive convolution of the CDFs of the $I_i$'s above with the CDFs of the integrals $\int_{t_{i-1}}^{t_i}\krnl(s) N(ds)$.   We know the smoothness of the CDF of each integral, but what about the smoothness of a convolution of two such CDFs?   The following corollary shows that the H\"{o}lder exponent associated to the convolution of any two of these CDFs will, at worst, be the lesser of the two H\"{o}lder exponents associated to the individual CDFs.  

\begin{cor}\label{c:sumIs}
Let $\krnl$, $n$ and $X(t)$ and $\cdf$ be as in Theorem (\ref{t:holder}), and suppose $X_0$ is a non-negative random variable independent of $X(t)$ with CDF $\cdf_0 \in C^{0,\gamma_0}(0,\infty)$.  Then, the CDF $\bar{\cdf}$ of the sum $X_0+ X(t)$ lies in the H\"{o}lder space $C^{0,\bar{\gamma}}(0,\infty)$, where $\bar{\gamma} = \min(\gamma,\gamma_0)$, and 
\begin{equation}
[\bar{\cdf}]_{C^{0,\bar{\gamma}}(0,\infty)} \leq \begin{cases} [\cdf]_{C^{0,\gamma}(0,\infty)} + \max(1,[\cdf_0]_{C^{0,\gamma_0}(0,\infty)}), \qquad \gamma_0 > \gamma \\
[\cdf]_{C^{0,\gamma}(0,\infty)}  + [\cdf_0]_{C^{0,\gamma_0}(0,\infty)},  \qquad \qquad \gamma_0 = \gamma \\
[\cdf_0]_{C^{0,\gamma_0}(0,\infty)} + \max(1,[\cdf]_{C^{0,\gamma}(0,\infty)}), \qquad \gamma_0 < \gamma \end{cases}
\end{equation}
\end{cor}

\begin{pf}
Let $0 < y < x$.  Observe first that if $\gamma > \bar{\gamma}$, then
\begin{eqnarray}
[\cdf]_{C^{0,\bar{\gamma}}} \leq \frac{\
\cdf(x) - \cdf(y)}{|x-y|^{\bar{\gamma}}} &\leq& \begin{cases} 1, \qquad |x - y| \geq 1 \\ \nonumber
\displaystyle \frac{\cdf(x) - \cdf(y)}{|x - y|^{\gamma}}, \quad |x-y| < 1 \end{cases} \\
 &\leq& \max(1, [\cdf]_{C^{0,\gamma}(0,\infty)} ) . \label{e:holderinq1}
\end{eqnarray}
Similarly, if $\gamma_0 > \bar{\gamma}$, we have
\begin{equation}
[\cdf_0]_{C^{0,\bar{\gamma}}} \leq \max(1, [\cdf_0]_{C^{0,\gamma_0}(0,\infty)} ) . \label{e:holderinq2}
\end{equation}

To compute the H\"{o}lder norm to  $\bar{\cdf}$, we write it as the convolution of $F$ and $F_0$:
\begin{equation}\label{e:fbar2p}
\bar{\cdf}(x) - \bar{\cdf}(y) = \int_0^\infty (\cdf(x - u) - \cdf(y-u)) dF_0(u).
\end{equation}
Since $\cdf = 0$ on the negative real axis, the above simplifies to
\begin{equation}
\bar{\cdf}(x) - \bar{\cdf}(y) = \int_0^y (\cdf(x - u) - \cdf(y-u)) dF_0(u) + \int_y^x \cdf(x - u) dF_0(u). \label{e:splitFbar}
\end{equation}
Now,  Theorem (\ref{t:holder}) and (\ref{e:holderinq1}) imply that
\begin{eqnarray}
\frac{|\int_0^y (\cdf(x - u) - \cdf(y-u)) d\cdf_0(u) |}{|x - y|^{\bar{\gamma}}} &\leq & \int_0^y \frac{|\cdf(x - u) - \cdf(y-u)|}{|x-y|^{\bar{\gamma}}} d\cdf_0(u)   \nonumber\\ 
& \leq & \begin{cases} [\cdf]_{C^{0,\gamma}(0,\infty)}, \qquad \gamma = \bar{\gamma} \\
\max(1, [\cdf]_{C^{0,\gamma}(0,\infty)} )\qquad \gamma > \bar{\gamma} \end{cases}. \label{e:gamgamp}
\end{eqnarray}
And, (\ref{e:holderinq2}) implies 
\begin{equation}
 \frac{\int_y^x \cdf(x - u) d\cdf_0(u)}{|x - y|^{\bar{\gamma}}} \leq \frac{\cdf_0(x) - \cdf_0(y)}{|x - y|^{\bar{\gamma}}} \leq \begin{cases}  [\cdf_0]_{C^{0,\gamma_0}(0,\infty)}, \qquad \gamma_0 = \bar{\gamma} \\
\max(1, [\cdf_0]_{C^{0,\gamma_0}(0,\infty)} )\qquad \gamma_0 > \bar{\gamma} \end{cases}. \label{e:gam0gamp}
\end{equation}
Thus, (\ref{e:splitFbar}), together with (\ref{e:gamgamp}) and (\ref{e:gam0gamp}) imply the result.
\end{pf}

Corollary \ref{c:sumIs} points to a possible ``worst case'' result about the
smoothness of a convolution of two CDFs with different smoothness
properties on $(0,\infty)$ by saying this convolution $\bar{F}$ on $\{x > 0\}$ will lie in the larger of
two H\"{o}lder spaces ($C^{0,\gamma}(0,\infty)$ and
$C^{0,\gamma_0}(0,\infty)$).  This worst case occurs if $\cdf$ and
$\cdf_0$ are both discontinuous at $x = 0$, as indicated in the
following corollary. 

\begin{cor}\label{c:worstcase}
Under the assumptions of Corollary \ref{c:sumIs}, suppose further that
$\cdf_0$ has a discontinuity at $x=0$ and that
$C^{0,\gamma}(0,\infty)$ and $C^{0,\gamma_0}(0,\infty)$ are the
smallest H\"{o}lder spaces to contain $\cdf$ and $\cdf_0$,
respectively.  Then,  $C^{0,\bar{\gamma}}(0,\infty)$ is the smallest
  H\"{o}lder space to contain $\bar{F}$.
\end{cor}

\begin{pf}

Notice that if $\cdf_0$ and $\cdf$ are discontinuous at $x=0$,  the quotients (\ref{e:gamgamp}) and
(\ref{e:gam0gamp}) can also be bounded below as
\begin{eqnarray*}
\frac{\int_0^y (\cdf(x - u) - \cdf(y-u)) d\cdf_0(u) }{|x -
  y|^{\bar{\gamma}}} &\geq& \frac{ \cdf(x) - \cdf(y)  } {|x -
  y|^{\bar{\gamma}}} (\cdf_0(0) - \cdf_0(0^{-})) \\
\frac{\int_y^x \cdf(x - u) d\cdf_0(u)}{|x - y|^{\bar{\gamma}}} &\geq&
\cdf(0)  \frac{\cdf_0(x) - \cdf_0(y)}{|x - y|^{\bar{\gamma}}} 
\end{eqnarray*}
And thus it follows from (\ref{e:splitFbar}) that if $C^{0,\gamma}$
and $C^{0,\gamma_0}$ are the smallest H\"{o}lder spaces to contain
$\cdf$ and $\cdf_0$, respectively,  then taking supremums over all $x$
and $y$ above will cause one of these two ratios to become infinite if
$\bar{\gamma} > \min(\gamma,\gamma_0)$.  Thus, (\ref{e:fbar2p})
implies that $C^{0,\bar{\gamma}}$ is the
smallest H\"{o}lder space to hold $\bar{F}$ on $\{x > 0\}$.   

\end{pf}

Notice that for sums like (\ref{e:poissummands}), the corresponding
CDFs will typically contain a jump at $x = 0$, thus the ``worst case''
is actually what occurs.

\section{Obtaining $\cdf(x,t)$}\label{s:methods}

There are standard finite-difference schemes for solving linear PDEs
(see for instance \cite{kincaid:1996}, Chapter 9).  
In this section, we describe such a scheme which can be used to solve (\ref{e:mainde})
iteratively on the interval $0 \leq t \leq T$.  We also study its convergence properties.

To begin, we need to define various components of a discrete approximation of $\cdf(x,t)$ which we will denote $\cdf(j \delta, i h)$ and define as follows.   Fix an interval $ [0,L] $ and subdivide it with a step size of $\delta$. That is, consider
\begin{eqnarray}
\Delta_{\delta} &=& \{  j \delta \ : \  0 \leq j \leq L/\delta \} \nonumber \\
&=& \{ x_0 < x_1 < x_2 \dots < x_M \},\label{e:defofDelta}
\end{eqnarray} 
where 
\begin{equation}\label{e:defofM}
M = |\Delta_\delta| = L/\delta +1
\end{equation}
 is the size of the mesh (assume $L/\delta$ is an integer).  We will also use a time step $h >0$ and consider the time points $\{ i h, \ i = 0,1,\dots,N \}$, where $N = T/h$.  Given a mesh $\Delta_\delta$, for each integer $ 0 \leq i \leq N $ define the column vector $\mathbf{\cdf}^i \in \mathbb{R}^M$ as 
\begin{equation}\label{e:colF}
\mathbf{\cdf}^i = \left(  \begin{array}{c}  
\cdf^i_0 \\
\cdf^i_1 \\
\cdf^i_2 \\
\vdots \\
\cdf^i_M \end{array} \right), 
\end{equation}
with $\cdf^i_j = \cdf(  j \delta, i h)$. In order to measure ``closeness'' on this mesh, we will use the following discrete $L^1$ norm defined on $\Delta_\delta$:
\begin{equation} \label{e:epsnorm}
\| \mathbf{u} \| _{1} = \sum_{j=0}^M \delta |u_j|, \qquad \mathbf{u} \in \mathbb{R}^M.
\end{equation}
Recall that for a $M \times M$ matrix $A = (a_{i,j})_{i,j=0}^M$,  the $L^{1}$ norm of $A$ is given by its maximal absolute column sum (\cite{kincaid:1996}, problem 4.4.11):
\begin{equation}\label{e:Anorm}
\| A \|_{1} = \max_{j=0,1,\dots,M} \sum_{i=0}^M |a_{i,j}|.
\end{equation}

\subsection{Finite-Difference Scheme for computing $\cdf(x,t)$}

\begin{figure}[ht]
\centering
\includegraphics[width=.8\textwidth]{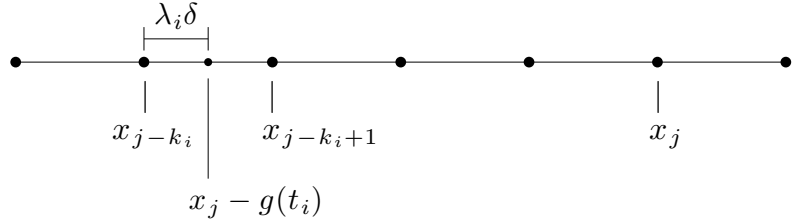}
\caption{ Relationship between $x_j$,  $\krnl(t_i)$, $k_i$, $\delta$ and $\lambda_i$. }
\label{f:FDcartoon}
\end{figure}

 We begin by rewriting equation (\ref{e:mainde}) as 
\begin{equation}\label{e:defoff}
{\cal P }[\cdf] \equiv \frac{\partial}{\partial t^+} \cdf(x,t) + n(t)(\cdf(x,t) - \cdf(x-\krnl(t),t)) = 0,
\end{equation}
where ${\cal P}$ is a linear operator defined on a suitable function space.  We must first choose a discrete approximation of ${\cal P}$.  For the time derivative, we use the usual forward difference approximation:
\begin{equation}
\frac{\partial}{\partial t^+}\cdf(x_j,t_i) \approx \frac{\cdf^{i+1}_j- \cdf^i_j}{h},
\end{equation}
where the $\cdf^{i}_j$ are defined in (\ref{e:colF}).  For the difference $\cdf(x,t) - \cdf(x-\krnl(t),t)$, we use a linear interpolation
\begin{equation}\label{e:linearinterp}
\cdf(x_j,t_i) - \cdf(x_j-\krnl(t_i),t_i) \approx \cdf^i_j - (1 - \lambda_i) \cdf^i_{j - k_i} + \lambda_i \cdf^i_{j - k_i + 1},
\end{equation}
where the integer $k_i$ satisfies $x_{j-k_i} < x_j - \krnl(t_i) < x_{j-k_i+1} $, and $\lambda_i  =   
\delta^{-1}( x_j - \krnl(t_i) -x_{j-k_i} )$.  In terms of $\delta$ and $h$, these are given by
\begin{equation}\label{e:defofk}
k_i = \left\lfloor \frac{\krnl(i h)}{\delta} \right\rfloor + 1,   \qquad \lambda_ i = \frac{|\delta k_i - \krnl(ih)|}{\delta}.
\end{equation}
See Figure \ref{f:FDcartoon}.

Putting these approximations together,  we define the \textit{forward difference} operator ${\cal P}_\krnl: \mathbb{R}^M \rightarrow \mathbb{R}^M$ as 
\begin{equation}\label{e:def_of_Ff}
({\cal P}_\krnl[\mathbf{\cdf}^i])_j = \frac{ \cdf^{i+1}_j - \cdf^{i}_j}{h} + n(i h) (\cdf^i_j - (1-\lambda_i) \cdf^i_{j - k_i} + \lambda_i) \cdf^i_{j - k_i + 1} )
\end{equation}
where $k_i,\lambda_i$ are defined in (\ref{e:defofk}).  In view of (\ref{e:defoff}), we require the sequence $\mathbf{\cdf}^i$ satisfy 
\begin{equation}\label{e:forwarditer}
{\cal P}_\krnl[\mathbf{\cdf}^i] = 0, \qquad i=1,2,\dots.
\end{equation}
From (\ref{e:def_of_Ff}), (\ref{e:forwarditer}) can we written as 
\begin{equation}\label{e:timestepup}
\cdf^{i+1}_j = \cdf^i_j - h n( i h )\cdf^i_j +  h  (1 - \lambda_i) \cdf^i_{j - k_i} + h \lambda_i  \cdf^i_{j - k_i + 1},
\end{equation}  
We can write this in matrix form.  Since $N$ is the size of time ($i$)
mesh and $M$ is the size of the space ($j$) mesh, we can express (\ref{e:timestepup}) as
\begin{equation}
\mathbf{\cdf}^{i+1} = A^{i} \mathbf{\cdf}^{i}, \qquad i=1,2,\dots,N,
\end{equation}
where $A^{i}$ is the $M\times M$ matrix\footnote{When implementing
  this method, $A^i$
  should be treated as a sparse matrix to avoid running out of
  memory} defined as
\begin{equation}\label{e:defofAmatrix}
A^i = (1 - n(i h) h) \mathbf{I}_{M \times M} +  \left( \begin{array}{cc} 
\mathbf{0}_{(k-1) \times (M-k+1)}  & \mathbf{0}_{(k-1) \times (k-1)} \\
\tilde{A}^i & \mathbf{0}_{(M-k+1) \times (k-1)} \end{array} \right),
\end{equation}  
where $\mathbf{I}_{M \times M}$ denotes the $M \times M$ identity matrix, and $\mathbf{0}_{M \times N}$ is a $M \times N$ zero matrix, and the $(M-k_i+1) \times (M-k_i+1)$ matrix $\tilde{A}^i$ is
\begin{equation}
(\tilde{A}^i)_{mn}  = \begin{cases}  h \lambda_i  n(i h), \quad m=n \\
h (1 - \lambda_i) n(i h) \quad m = n+1 \\
0 \quad \mathrm{else} \end{cases}.
\end{equation}
With this, $\mathbf{\cdf}^i$ is easily calculated as
\begin{equation}
\mathbf{\cdf}^i = A^{i} A^{i-1} \dots A^1 \mathbf{\cdf}^0.
\end{equation}
where the initial condition $\mathbf{\cdf}^0$ is determined by the CDF of $X_0$ and is known exactly.

\subsection{Error Analysis}\label{s:error}

In this section, we study the error associated with solving
differential-difference equation (\ref{e:mainde}) with the
finite-difference method described in the previous section.  We'll
require the following assumptions on the kernel $\krnl$ and control measure $n$:
\begin{align}
&  \mbox{$\bullet$   $\krnl$   is a monotone increasing function with H\"{o}lder
  continuous inverse $\krnl^{-1}$ } \notag \\
& \mbox{$\bullet$  $\krnl$ and $n$ are
  Lipschitz functions with} \label{e:gnassump} \\
&  \ \ \ |\krnl(t) - \krnl(s)| \leq L_\krnl |s - t|,
  \quad |n(s) - n(t)| \leq L_{n} |s - t| \    \mbox{for constants $L_\krnl,L_n \geq 0$} \notag.
\end{align}

We will see that our
method gives accuracy of order $o( \delta^\gamma + h^{\gamma})$ in the
discrete $L^1$ norm,  where $\gamma$ is the H\"{o}lder exponent of
$\krnl^{-1}$. As above, $\delta$ is the size of the spacial mesh, and
$h$ is the size of the temporal mesh.  

We first consider the error associated with the approximations used in our discretization:
\begin{eqnarray*}
\frac{\partial \cdf}{\partial t^{+}} (x_j,t_\ell) &\approx& \frac{\cdf(x_j,t_\ell+h) - \cdf(x_j,t_\ell)}{h} \\
\cdf(x_j - \krnl(t_\ell),t_\ell) &\approx&  (1 - \lambda_\ell) \cdf(x_{j - k_\ell},t_\ell) + \lambda_\ell \cdf(x_{j - k_\ell+1},t_\ell).
\end{eqnarray*}
Consider the meshes $0 = x_0 < x_1 < x_2 < \dots x_M$ and $0 = t_0 < t_1 < \dots < t_N$. For integers $0 \leq j \leq M$ and $0 \leq \ell \leq N$, define the absolute differences
\begin{eqnarray}
S_{j,\ell}(h) &=& \cdf(x_j,t_\ell + h) - \cdf(x_j,t_\ell) - h \frac{\partial \cdf}{\partial t^{+}}(x_j,t_\ell)   \label{e:defofS}  \\
R_{j,\ell}(\delta) &=&    ( (1 - \lambda_\ell) \cdf(x_{j - k_\ell},t_\ell) + \lambda_\ell \cdf(x_{j - k_\ell + 1} , t_\ell ))  - \cdf(x_j - \krnl(t_\ell),t_\ell)    . \label{e:defofR}
\end{eqnarray}

In order to ensure convergence in our finite-difference scheme, these two errors must approach $0$ ``fast enough'' as $\delta,h \rightarrow 0$.  Bounds for $S_{j,\ell}$ and $R_{j,\ell}$ are given in Lemmas \ref{l:Slem} and \ref{l:Rlem} in section \ref{s:approx}.   In short, these Lemmas imply the following:
\begin{eqnarray*}
|S_{j,\ell}(h)| &=& \begin{cases}  o(h^{1+\gamma}) \quad \mbox{if } x_j - g(t_{\ell+1}) \geq 0 \\
o(h) \quad \mbox{if } x_j - g(t_{\ell}) \geq 0 > x_j - g(t_{\ell + 1}) \\
o(h^2) \quad \mbox{if } x_j - g(t_\ell) < 0 \end{cases},  \\
|R_{j,\ell}(\delta)| &=& o(\delta^{\gamma}).
\end{eqnarray*}

For $\ell \geq 0$, let $\mathbf{\cdf}_e^\ell \in \mathbb{R}^M$ be the vector of {\it exact} values of $\cdf$ at time $t_\ell$, i.e.
\begin{equation}\label{e:Pexact}
 \mathbf{\cdf}_e^\ell = \left( \begin{array}{c}
\cdf(x_0,t_\ell) \\
\cdf(x_1,t_\ell) \\
\vdots \\
\cdf(x_M,t_\ell) \end{array} \right).
\end{equation}
We will now state a theorem which shows that our method converges to
the exact solution to (\ref{e:mainde}) in the discrete $L^1$ norm with
a convergence rate faster than a constant multiple of $(h^\gamma + \delta^\gamma)$, i.e.~that
$\| \mathbf{\cdf}^\ell - \mathbf{\cdf}^\ell_e\|_1 = o(h^\gamma +
\delta^\gamma)$ and $\delta, h \rightarrow 0$.

As the proof of the following theorem will show, we will require the following {\it stability criterion} on $h$ to guarantee convergence:
\begin{equation*}
\fbox{ $\mbox{Stability Criterion:      } \ \   h n^{\star}_T < 1.$ }
\end{equation*}
where $\displaystyle n^{\star}_T = \sup_{0 \leq s \leq T} n(s)$ and
$h$ is the time step used in the finite-difference scheme.  This
can always be met easily if $n$ is bounded.  This condition makes
sense since the differential equation (\ref{e:mainde}) holds because the probability of 2 or more Poisson arrivals in an infinitesimal time
interval is negligible.  Hence, when approximating (\ref{e:mainde}) with
a finite-difference method, we
must ensure that the time step we consider is small enough relative
to the control measure to make this probability small.

\begin{thm}\label{t:errorbound}
Let $\krnl$ be a positive, continuous, and strictly increasing kernel with inverse $\krnl^{-1} \in C^{0,\gamma}(\krnl(0),\krnl(T))$ for some $0 < \gamma \leq 1$.  Suppose the control measure $n(s)$ is a bounded function and set $n^{\star}_T = \sup_{0 \leq t \leq T} n(t)$.  Let $\mathbf{\cdf}_e^i$ be defined as in (\ref{e:Pexact}) and let $\mathbf{\cdf}^i$ $i=1,2,\dots,N$ be the approximations given by the forward difference scheme.  Then, if $h n^{\star}_T < 1$, the forward difference scheme is convergent and satisfies 
\begin{equation}\label{e:errorbound}
\| \mathbf{\cdf}^N - \mathbf{\cdf}_e^N \|_{1} = o(h^\gamma + \delta^\gamma)
\end{equation}

 \end{thm}

\begin{pf}

At each step in the iteration, write $A^i \mathbf{\cdf}^i_e = \mathbf{\cdf}^{i+1}_e + \mathbf{\epsilon}_i$, where $\epsilon_i \in \mathbb{R}^n$ is the error introduced at step $i$.  Then, since the initial condition $\mathbf{\cdf}^0$ is known exactly, we have
\begin{eqnarray}
\mathbf{\cdf}^0 &=& \mathbf{\cdf}^0_e  \nonumber \\
\mathbf{\cdf}^1 &=& A^0 \mathbf{\cdf}^0 = \mathbf{\cdf}^1_e + \epsilon_0  \nonumber\\
\mathbf{\cdf}^2 &=& A^1 \mathbf{\cdf}^1 = \mathbf{\cdf}^2_e + \epsilon_1+ A^1 \epsilon_0 \nonumber\\
&\vdots& \nonumber\\
\mathbf{\cdf}^N &=& \mathbf{\cdf}^N_e + \epsilon_{N-1} + A^{N-1} \epsilon_{N-2} + \dots + A^{N-1}A^{N-2}\dots A^{1} \epsilon_0 \label{e:errorseries}
\end{eqnarray}
After subtracting $\mathbf{\cdf}^e$ from both sides and taking norms, (\ref{e:errorseries}) implies
\begin{eqnarray}
\| \mathbf{\cdf}^N - \mathbf{\cdf}_e^N \|_{1} &\leq& \sum_{\ell=0}^{N-1} \left\| \prod_{i=\ell+1}^{N-1} A^{i}  \right\|_{1}  \| \epsilon_\ell \|_1  \nonumber\\
&\leq& \sum_{\ell=0}^{N-1} \left( \prod_{i=\ell+1}^{N-1} \left\|  A^{i}  \right\|_{1}  \right)\| \epsilon_\ell \|_1 \label{e:boundstep1}
\end{eqnarray}
where $\| \epsilon_\ell \|_1$ and $\left\|  A^{i}  \right\|_{1}$ are
defined respectively in (\ref{e:epsnorm}) and (\ref{e:Anorm}).
From the definition of $A^i$ in (\ref{e:defofAmatrix}) and the stability condition $n^{\star}_T h < 1$, the norm $\|A^i\|_{1}$ is given by the maximum absolute sum of the columns of $A^{i}$, which in our case is
\begin{equation}
\|A^i\|_{1} = |1 - n(t_i) h| + |h n(t_i) \lambda_i | + |h n(t_i) (1- \lambda_i) | = 1.
\end{equation}
Thus, (\ref{e:boundstep1}) becomes
\begin{equation}\label{e:boundstep2}
\| \mathbf{\cdf}^i - \mathbf{\cdf}_e^i \|_{1}  \leq \sum_{\ell=0}^{N-1} \| \epsilon_\ell \|_1.
\end{equation}

What remains is to bound 
\begin{equation*}
\| \epsilon_\ell \|_1 = \| \mathbf{\cdf}^{\ell+1}_e- A^{\ell} \mathbf{\cdf}_e^\ell  \|_1 = \sum_{j=0}^M \delta |(\epsilon_\ell)_j|
\end{equation*}
Using the differential equation (\ref{e:mainde}) and the definitions
of the errors $S_{j,\ell}$ and $R_{j,\ell}$ in (\ref{e:defofS}) and
(\ref{e:defofR}), we have
\begin{eqnarray}
(\epsilon_\ell)_j &=& \cdf(x_j,t_{\ell + 1}) - (A^\ell \mathbf{\cdf}_e^\ell )_j \nonumber \\
&=&  \cdf(x_j,t_{\ell + 1}) - \{ \cdf( x_j , t_\ell ) ( 1 - h
n(t_\ell) ) +  h n(t_\ell) [ (1- \lambda_\ell) \cdf(x_{j - k_\ell},t_\ell) + \lambda_\ell \cdf(x_{j - k_\ell+1},t_\ell) ] \} \nonumber \\
& = &  \cdf(x_j,t_\ell + h) - \cdf(x_j,t_\ell)  - h [ n( t_\ell ) (
(1- \lambda_\ell) \cdf(x_{j - k_\ell},t_\ell) + \lambda_\ell \cdf(x_{j
  - k_\ell+1},t_\ell)  - \cdf(x_j,t_\ell)) ] \nonumber \\
&=&  \cdf(x_j,t_\ell + h) - \cdf(x_j,t_\ell)  - h [ n( t_\ell ) (
 \cdf(x_j- \krnl(t_\ell),t_\ell) 
 - \cdf(x_j,t_\ell) +R_{j,\ell}(\delta)  ) ] \nonumber \\
&=&   \cdf(x_j,t_\ell + h) - \cdf(x_j,t_\ell)  - h [ \frac{ \partial \cdf}{\partial t^+}(x_j,t_\ell) + n(t_\ell) R_{j,\ell}(\delta) ]  \nonumber \\
&=& S_{j,\ell}(h) -  n(t_\ell) R_{j,\ell}(\delta) h .     \label{e:epsilonbound1}
\end{eqnarray}

Now, for $\ell$ fixed, consider the following partition of the mesh
$\Delta_\delta$ in (\ref{e:defofDelta}):
\begin{eqnarray*}
I_1^{(\ell)} &=& \{ j  \ : \  x_j - \krnl(t_\ell) \geq 0 > x_j - \krnl(t_{\ell + 1})  \} \\
I_2^{(\ell)} &=& \{ j \ : \ x_j - \krnl(t_\ell) < 0 \} \\
I_3^{(\ell)} &=& \{ j \ : \ x_j - \krnl(t_{\ell+1}) \geq 0 \}
\end{eqnarray*}
Notice that $I_1^{(\ell)}$ contains at most one element.  Using
(\ref{e:epsilonbound1})  combined with Lemmas \ref{l:Slem} and
\ref{l:Rlem} and that $M \delta = L$, we have
\begin{eqnarray}
\| \epsilon_\ell \|_1 = \sum_{j=0}^M\delta |(\epsilon_\ell)_j|  &\leq&  n(t_\ell) h \sum_{j=0}^M \delta |R_{j,\ell}(\delta) | + \sum_{j \in I_1^{(\ell)}} \delta |S_{j,\ell}(h)| + \sum_{j \in I_2^{(\ell)}} \delta |S_{j,\ell}(h)| +  \sum_{j \in I_3^{(\ell)}} \delta |S_{j,\ell}(h)| \nonumber \\
&\leq&   O( h \delta^{\gamma}) + | I_1^{(\ell)} | \cdot  O( \delta h ) + | I_2^{(\ell)} | \cdot O(\delta h^2) + | I_3^{(\ell)} | \cdot O(\delta h^{1+\gamma}) \nonumber  \\
&=& O(h \delta^\gamma) + O(\delta h) + O(h^2) + O(h^{1+\gamma}) \nonumber \\
&=& O(h \delta^\gamma ) + O(h^{1+\gamma}) . \label{e:bigOs}
\end{eqnarray}
Finally, since $N h = T$, (\ref{e:boundstep2}) and (\ref{e:bigOs}) imply
\begin{eqnarray}
\| \mathbf{\cdf}^i - \mathbf{\cdf}_e^i \|_{1}  &\leq& \sum_{\ell=0}^{N-1} \| \epsilon_\ell \|_1. \nonumber \\
&\leq & O(h^\gamma + \delta^\gamma),  
\end{eqnarray}  
which finishes the proof.

\end{pf}

\subsection{Extension to more general kernels}\label{s:extend}

In order to guarantee the error bound (\ref{e:errorbound}), we must
assume that $\krnl$ is a strictly  increasing, continuous function
whose inverse is H\"{o}lder continuous.  This assumption is somewhat
restrictive, as we would like to apply our method to a wider set of
kernels.   In this section, we indicate how to transform a  problem with a general kernel into problem which satisfies the assumptions of Theorem \ref{t:errorbound}.  

Recall that if $X,Y$ are  independent non-negative random variables with respective CDFs $\cdf_X, \cdf_Y$, then the CDF of the sum $X+Y$ is given by the convolution
\begin{equation}\label{e:convolve}
\cdf_{X+Y}( u ) = \int_0^\infty \cdf_Y( u - x ) \cdf_X(dx).
\end{equation}
If $\cdf_X$ and $\cdf_Y$ are both defined on a mesh $ \Delta_{\delta} = \{x_j\}_{j=1}^M$, then the convolution (\ref{e:convolve}) can be approximated by the discrete convolution
\begin{equation}\label{e:discreteconvolve}
\cdf_{X+Y}(u) \approx \sum_{j= 0 }^M  \cdf_Y( u - x_j ) [ \cdf_X( x_{j+1} ) - \cdf_X(x_j) ].
\end{equation}

\subsubsection*{Decreasing kernels}
If $\krnl$ is a strictly {\it decreasing} kernel with H\"{o}lder continuous inverse, then a simple variable transform re-expresses the integral into one with a strictly increasing kernel:
\begin{equation}
I(\krnl) = \int_0^T\krnl(s)N(ds) \overset{d}{=} \int_0^T \widetilde{f}(s) \widetilde{N}(ds)
\end{equation}
where $\widetilde{f}(s) = \krnl(T-s)$ and $\widetilde{N}(ds)$ is a Poisson random measure with control $\widetilde{n}(s) = n(T-s)$.  

\subsubsection*{Flat kernels}
If $\krnl \equiv \lambda$ is a constant function, then there is no need to approximate the CDF since its distribution is known exactly:
\begin{equation}
I(\krnl) = \int_0^T \krnl N(ds) \sim  \lambda  \mathrm{Pois}( \int_0^T n(s) ds  ).
\end{equation}

\subsubsection*{Negative kernels}
If $\krnl$ is a negative strictly monotone increasing/decreasing kernel with H\"{o}lder continuous inverse, then the negative of the integral fits the proper assumptions:
\begin{equation}\label{e:convolveP}
I(\krnl) = \int_0^T \krnl(s) N(ds) \overset{d}{=} - \left( \int_0^T (-\krnl(s)) N(ds) \right).
\end{equation}
Thus, ones computes the CDF of $-I(\krnl)$ with the forward difference
method, and then uses the relationship $\cdf_{I(\krnl)}(u) =1 -  \cdf_{-I(\krnl)}(-u)$.  

\subsubsection*{``Piecewise'' kernels}
Combining the above three ideas, the convolution formula (\ref{e:discreteconvolve}) and the independent increment property of the random measure $N$, we can approximate integrals with {\it piecewise} increasing/flat/decreasing, positive/negative kernels, whose inverse on each non-flat piece is H\"{o}lder continuous.  Indeed, if $\krnl$ is  increasing/flat/decreasing and positive/negative between each of the points $0 = t_0 < t_1 < \dots < t_{n-1} < t_n = T$, then we have
\begin{equation}
I(\krnl)  = \int_0^T \krnl(s) N(s) \overset{d}{=} \sum_{i=0}^{n-1} \int_{t_i}^{t_{i+1}} \krnl(s) N(ds) 
\end{equation}
This sum on the right-hand side can be computed inductively by using the integral up to $t_{i-1}$ as an initial condition, and iterating until the CDF is computed up to $t_i$.  That is, set $I_0 = 0$, and use our method to compute
\begin{equation}
I_i = I_{i-1} + \int_{t_{i-1}}^{t_i} \krnl(s) N(ds), \quad i=1,2,\dots,n.
\end{equation}
Corollary \ref{c:sumIs} and Theorem \ref{t:errorbound} imply that the rate of convergence in this
method in the $L^1$ norm will depend on, at worse, the minimum of the H\"{o}lder exponents of $\krnl^{-1}$ over each of the intervals $[t_i,t_{i+1})$.

\section{Examples}\label{s:examples}

We shall demonstrate our method on some examples.   We start with 
\begin{equation}\label{e:simpex}
\int_0^1 s N(ds),
\end{equation}
where $N$ is a Poisson random measure with control Lebesgue: i.e.~ $n(s) = 1$.  This example is simple enough to compute exactly, so it serves as an useful test case.  The second example is where the kernel $\krnl$ is a parabola
\begin{equation}
\int_0^2 (1 - (1 - s)^2) N(ds),
\end{equation}
 and the control measure is again Lebesgue.  In this example, the
 kernel $\krnl$ has an inverse which lies in the H\"{o}lder space
 $C^{0,1/2}(0,1)$, and thus we expect the CDF to share this property.
 For purposes of comparison, we also compute this CDF by approximating the integral
 (\ref{e:abatesmet}).  

  Finally, we consider an
 integrand which is both positive and negative:
\begin{equation}
\int_0^1 \sin(2 \pi s) N(ds)
\end{equation}
with control measure Lebesgue.  This will demonstrate some of the
ideas in Section \ref{s:extend}.

Recall that in each example, we will see a jump at $x=0$ of size
$\exp(-\int_0^T n(s) ds)$ in the CDF which is the probability of no
Poisson arrivals.

\subsubsection*{Example 1:   $\displaystyle   \int_0^1 s N(ds), \quad n(s) = 1 $}

In this example, the kernel is the identity function $\krnl(s) = s$,
hence the integral is given by the sum of arrival times in $(0,1)$.
The number of arrivals follows a Poisson distribution with rate
$\int_0^1 n(s) ds = 1$
and in the event of $k \geq 0$ arrivals, the times of these arrivals are given by $k$
i.i.d.~random variables $U_1,U_2,\dots,U_k$ which are uniform on
$(0,1)$.  Thus,  the CDF $\cdf$ of $\int_0^1 s N(ds)$ can be expressed by conditioning on the number of Poisson arrivals in the interval $(0,1)$:
\begin{equation}\label{e:exactform1}
\cdf(x) = \sum_{k=0}^\infty \mathrm{Prob} \left. \left[ \sum_{i=1}^k U_i \leq
  x \right| \mbox{$k$ arrivals} \right] \mathrm{Prob}[ \mbox{$k$ arrivals}]=   \sum_{k=0}^\infty \frac{P_k(x) e^{-1}}{ k! },
\end{equation}
where $P_k(x)$ is the CDF of a sum of $k$ i.i.d.~$U(0,1)$ random variables.  Each $P_k$ can written as a piecewise polynomial, and can be computed recursively as
\begin{eqnarray*}
P_0(x) &=& \begin{cases} 0 \quad x < 0 \\
1 \quad x \geq 0 \end{cases}  \\
P_k(x) &=& \int_0^1 P_{k-1}(x-u) du \quad k \geq 1.
\end{eqnarray*}
The first few $P_k$ are
\begin{equation*}
P_1(x) = \begin{cases} 0, \quad x < 0 \\
x, \quad 0 \leq x < 1 \\
1, \quad x \geq 1 \end{cases}, 
\quad P_2(x) = \begin{cases} 0, \quad x < 0 \\
\frac{x^2}{2}, \quad 0 \leq x < 1 \\
\frac{(-2 + 4 x - x^2)}{2}, \quad 1 \leq x < 2 \\
1, \quad x \geq 2
\end{cases}, \quad
P_3(x) = \begin{cases}
0, \quad x < 0 \\
\frac{x^3}{6}, \quad 0 \leq x < 1 \\
 \frac{ 3 - 9 x + 9 x^2 - 2 x^3}{6}, \quad  1 \leq x < 2 \\
 \frac{-21 + 27 x - 9 x^2 + x^3}{6}, \quad  2 \leq x < 3 \\
1, \quad x \geq 3 \end{cases}.
\end{equation*}
Notice that it suffices to take the first $11$ terms in the sum (\ref{e:exactform1}), since $\sum_{i=11}^\infty e^{-1}/n! < 10^{-7}$, which is much more than the precision we seek.  

\begin{figure}[h]
\centering
\includegraphics[scale=.5]{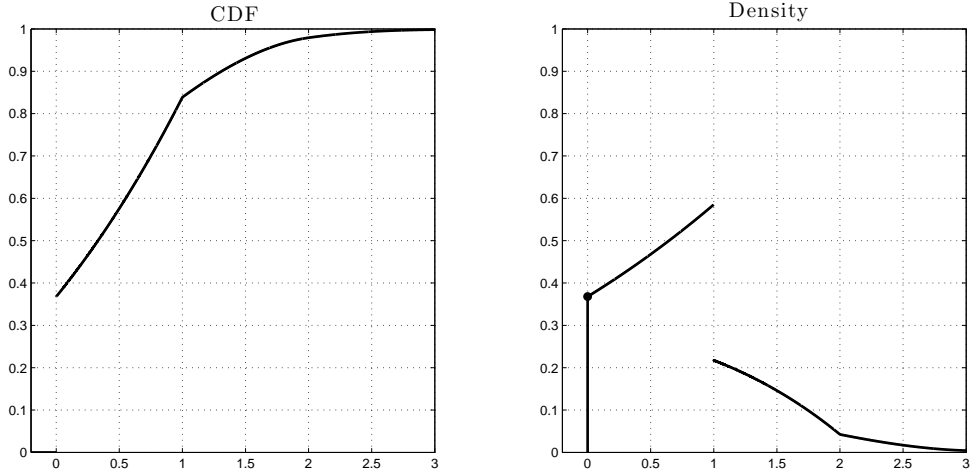}
\caption{CDF and density for the stochastic integral (\ref{e:simpex}).
  The vertical line with a dot in the plot of the density represents a
  discrete atom in the distribution whose mass is the height of the vertical line.  To compute this, we used the forward difference method with $h = \delta =  5 \times 10^{-4}$.  On a 1.5 Ghz computer, this computation took less than 5 seconds using MATLAB. Notice the ``kink'' in the CDF and discontinuity in the density at $x=1$.  This is not surprising, since here $\gamma = 1$ so expect the CDF to be at worst Lipschitz continuous. } \label{f:intsNds}
\end{figure}

We generated $\cdf(x)$ on the interval $ x \in [0,3)$ with our method with
$\delta = h =10^{-4}$.  This is plotted in Figure \ref{f:intsNds}
along with a numerical approximation of the density\footnote{This is
  obtained using the central difference $\cdf'(x) \approx
  (2 \delta_1)^{-1} ( F(x + \delta_1) - F(x-\delta_1))$.  We must choose
  $\delta_1 > \delta$ to account for the error made in computing the
  CDF with step size $\delta$.}.
Notice the ``kink'' at $x=1$ in the CDF and the corresponding
discontinuity in the density.  This is consistent with Theorem
\ref{t:holder},  since here the inverse of the kernel is simply
$\krnl^{-1}(s) = s$, which lies in the space $C^{0,1}(0,1)$, so here
we take $\gamma = 1$.  Thus, we
expect that the $\cdf$ is at worst Lipschitz continuous.   This is
indeed the case, since $\cdf$ is a linear combination of the $P_n$'s,
and $P_1$ is clearly in $C^{0,1}(0,\infty)$ and no smaller H\"{o}lder
space.  

\begin{figure}[h]
\centering
\includegraphics[scale=.5]{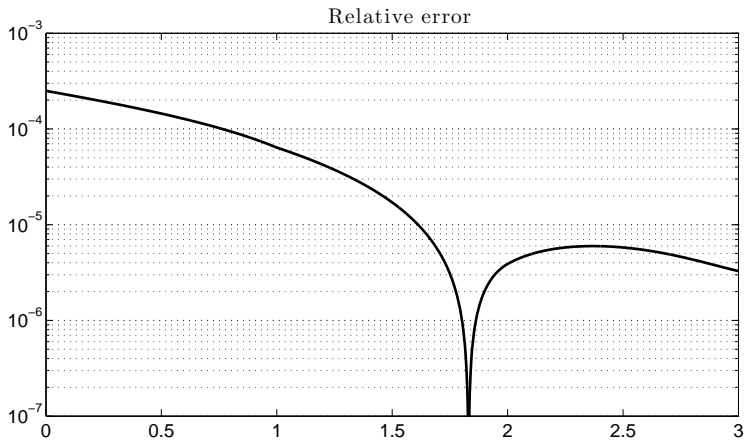}
\caption{ The relative error between the CDF of (\ref{e:simpex}) computed using our method, and the true value of the CDF computed with (\ref{e:exactform1}).  We see that our method gives accuracy of $\ll .1 \%$ over all $x$ considered.    } \label{f:intsNdserr}
\end{figure}

Using the first $11$ terms in the sum (\ref{e:exactform1}), we also computed this CDF on the interval $[0,3)$ and looked at the relative error between this and the result of our method.   The relative error is plotted in Figure \ref{f:intsNdserr}.  We have excellent agreement over all $x$ with our method.

\subsubsection*{Example 2: $\displaystyle  \int_0^2 [1 - (1-s)^2] N(ds)$, $\quad n(s) = 1$}

\begin{figure}[h]
\centering
\includegraphics[width=1.0\textwidth]{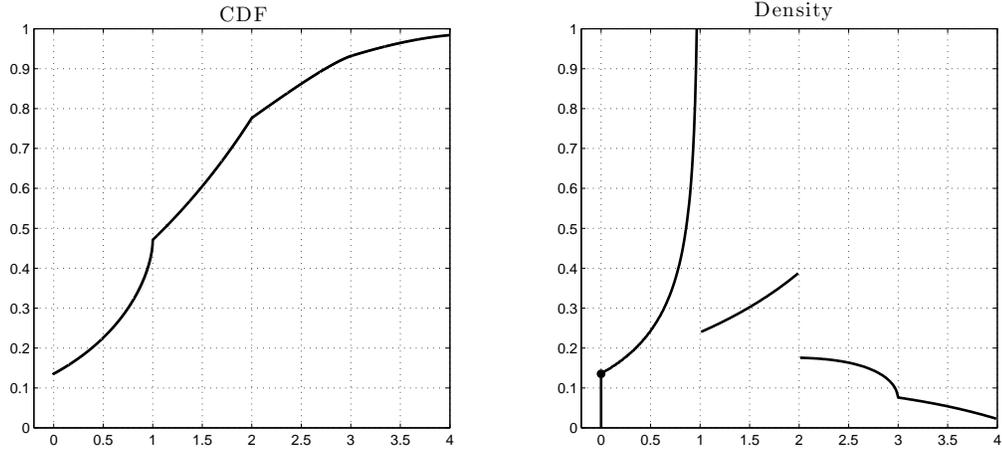}
\caption{ CDF and density the Poisson integral of $g$.  In this case, the CDF is H\"{o}lder continuous for $x > 0$ with H\"{o}lder exponent $\gamma = 1/2$.  As a result, we see a sharp ``corner'' in the CDF and a singularity in the density at $x = 1$.  To generate the CDF, we used $\delta = h = 10^{-4}$.  }
\label{f:intparNds}
\end{figure}

Here we consider the more complicated example where $\krnl(s) = 1 -
(1-s)^2$ for $0 \leq s \leq 2$, which is a downward facing parabola with maximum value $1$ and roots at $0$ and $2$.   This function is piecewise monotone on the intervals $(0,1)$, and $(1,2)$, so the finite difference scheme can be applied first to compute $F(x,t)$ for $0 \leq t \leq 1$, and then subsequently for $1 \leq t \leq 2$ by using $F(x,1)$ as an initial condition.  The density can be obtained with numerical differentiation.  This is shown in Figure \ref{f:intparNds}.

Notice the sharp corner in the CDF at $x = 1$ and the corresponding
singularity in the density.  This illustrates the results of Theorem
\ref{t:holder} and Corollary \ref{c:worstcase}, that is, the CDF lies in the smallest H\"{o}lder space $C^{0,\gamma}(0,\infty)$ which  contains the inverse of $g$ restricted to $(0,1)$ and the inverse of $g$ restricted to $(1,2)$.  In this case we have
\begin{eqnarray*}
g^{-1}(x) &=& 1 - \sqrt{1-x},  \quad   \mbox{for $g$ restricted to $(0,1)$}  \\
g^{-1}(x) &=& 1 + \sqrt{1-x},  \quad  \mbox{for $g$ restricted to $(1,2)$} 
\end{eqnarray*}
In both of these cases, $g^{-1} \in C^{0,1/2}(0,1)$,  and thus we
expect the CDF $F(x,t) \in C^{0,1/2}(0,\infty)$ as well.

For comparison purposes, we also computed this CDF using the
integration method of Abate and Whitt discussed briefly in the
introduction.  Since this integral is only applicable for continuous $\cdf$, we first rewrite the characteristic
function $\phi_g$ by ``removing'' the case where there are 0 arrivals,
and considering this case separately,
\begin{eqnarray*}
\phi_g(\theta) &=& \mathbb{E} \exp\left( i \theta \int_0^2 \krnl(s) N(ds)   \right) \\
&=& \mbox{Prob}[N=0] e^{0} + \left. \mathbb{E} \left[ \exp\left(  i \theta
    \int_0^2 \krnl(s) N(ds)   \right) \right|  N \geq 1 \right] \\
&=& e^{-2}  + \left[ \exp\left( \int_0^2 (e^{i \theta \krnl(s)}
  - 1) ds\right) - e^{-2} \right] \\
&\equiv& e^{-2} + \bar{\phi}_g(\theta).
\end{eqnarray*}
Then, using (\ref{e:abatesmet}) the CDF $\cdf$ is given by
\begin{equation}\label{e:phibar}
\cdf(x) = e^{-2} + \int_0^\infty \mbox{Re}( \bar{\phi}_g )(u) \frac{
  \sin( x u) }{u} du.
\end{equation} 
To approximate the integral in (\ref{e:phibar}), we truncate the
infinite limit at $T_I >0$, choose a mesh size $\eta> 0$
and use a trapezoid approximation:
\begin{equation}\label{e:trap}
\cdf(x) \approx e^{-2} + \frac{\eta}{2} \left( \mbox{Re}(\bar{\phi}_g )(0) +
\mbox{Re}(\bar{\phi}_g)(T_I) \frac{\sin(T_I x)}{T_I} \right)  + \eta
\sum_{j=1}^{K-1}   \mbox{Re}(\bar{\phi}_g)(j \eta) \frac{\sin(j \eta
  x)}{ j \eta } , 
\end{equation}
Where $K = T_I/\eta$ is the size of the mesh.  Notice that computing $
\mbox{Re}(\bar{\phi}_g)(j \eta)$ for every $j$ involves also computing
the integral $\int_0^2 (e^{i \theta \krnl(s)} - 1) ds $.  To do so, we used the MATLAB function
  \texttt{quad}.  The integrand is a highly oscillating function for large
$\theta$, causing numerical approximation methods to converge slowly.  This makes
computing the sum in (\ref{e:trap}) {\it extremely} time-consuming for large
$K$.

\begin{figure}[h]
\centering
\includegraphics[width=1.0\textwidth]{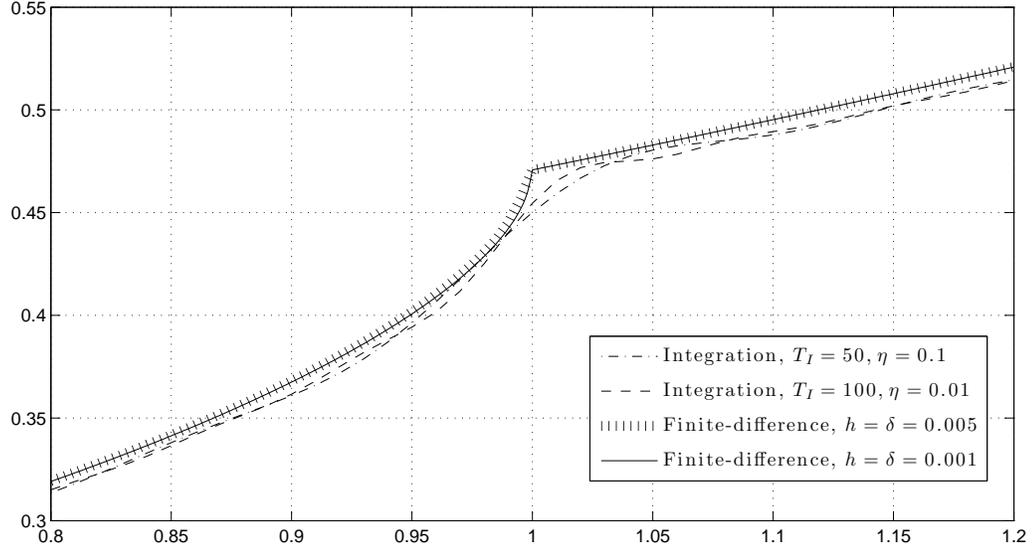}
\caption{ A comparison between the CDF of Example 2 computed using the
integration method (\ref{e:trap}) and the finite-difference scheme of
Section \ref{s:methods}, zoomed in on the kink in the CDF at $x=1$.
The integration method takes over 200 times 
longer to implement, and it seems that even larger $T_I$ and smaller $\eta$
is required to obtain sufficient convergence.  Also, important
features of $\cdf$ are missed with integration, such as the kink at $x=1$ which is
``smoothed'' here.  The finite-difference
method, on the other hand, converges quickly in $\delta$ and $h$, 
is much faster, and captures the kink in $\cdf$}
\label{f:FDintcompare}
\end{figure}

In Figure \ref{f:FDintcompare}, we plotted the CDF around the kink at $x=1$ using these two
methods.  We used the finite-difference method with $\delta = h =
 5 \times 10^{-3}$ and $\delta = h =1 \times 10^{-3}$, and used the integration method
with $T_I = 50$, $\eta = .1$ and $T_I = 100$ and $\eta = .01$.  Notice
that there is a small difference between the curves generated by
the finite-difference method, which implies adequate convergence has
been met.   On
the other hand, the results of the integration have ``smoothed'' out
the kink, implying a larger $T_I$ is required.  Moreover, the computation time
for these methods varied drastically, to generate the curves in Figure
\ref{f:FDintcompare}, the finite-difference method required less than 2 seconds, while
integration requires about 8 minutes.

\subsubsection*{Example 3: $\displaystyle  \int_0^1 \sin(2 \pi s) N(ds)$, $\quad n(s) = 1 $}

\begin{figure}[h]
\centering
\includegraphics[width=1.0\textwidth]{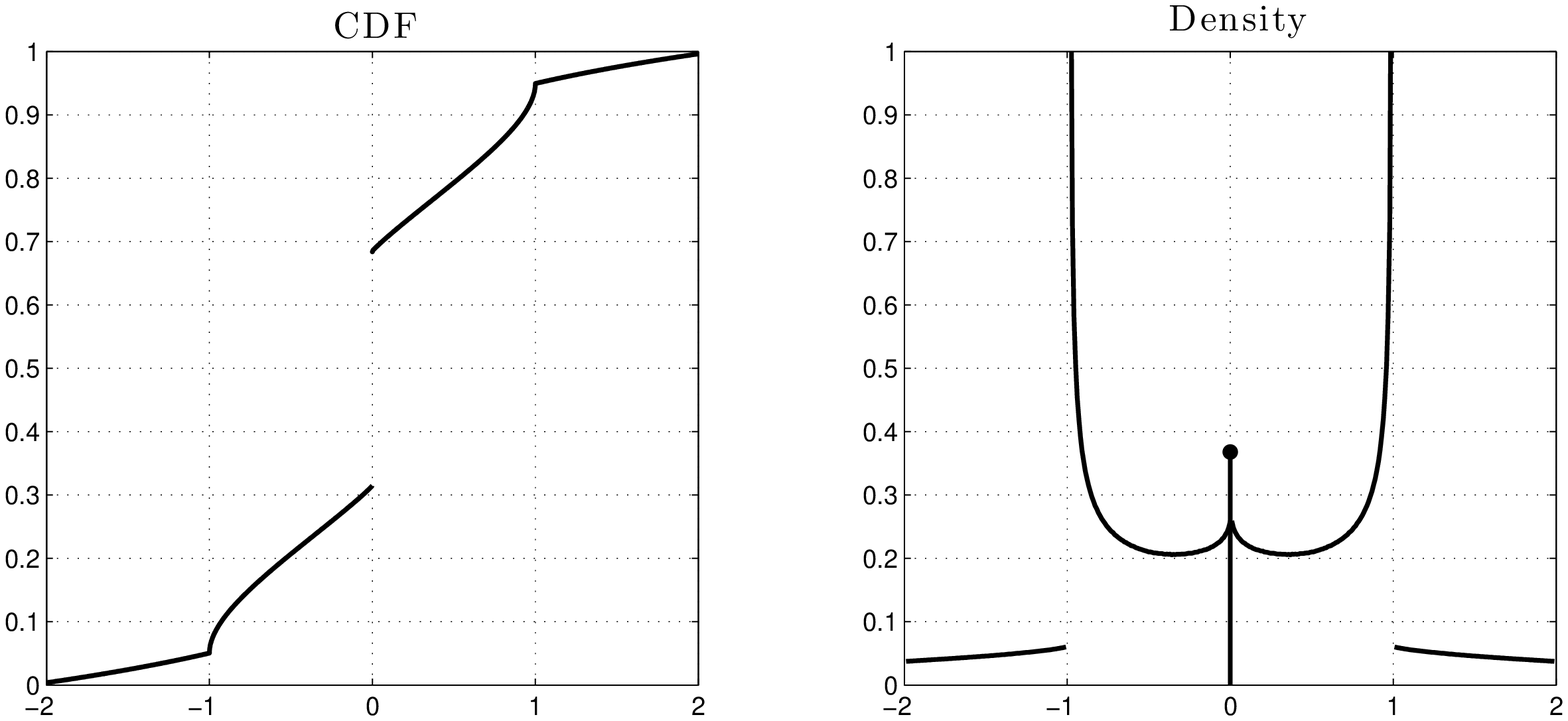}
\caption{ CDF and density the Poisson integral of the sine wave $\krnl(s) = \sin(2 \pi s)$ over the interval $0 \leq s \leq 1$. }
\label{f:intsinsNds}
\end{figure}

We now consider the Poisson integral of one period of a sine wave which is a case where the kernel is both positive and negative.  To compute the CDF $F$, we split up the integral as
\begin{equation}\label{e:breakitupP}
\int_0^1 \sin(2 \pi s) N(ds) = \int_0^{\frac{1}{2}}   \sin(2 \pi s) N(ds) + \int_{\frac{1}{2}}^1   \sin(2 \pi s) N(ds) =: I_1 + I_2
\end{equation} 
and note that the summands on the right are independent and satisfy $I_1 \overset{d}{=} -I_2$.  Thus, if $F_1(x)$ is the CDF of $I_1$, then $F_2(x) = 1 - F_1(-x)$ gives the CDF of $I_2$.  Taking the convolution of $F_1$ and $F_2$ we use the symmetry of $I_1$ and $I_2$ and see
\begin{eqnarray}
F(u) &=& \int_{-\infty}^\infty F_1(u - x) dF_2(x) \nonumber \\
&=& e^{-1/2} F_1(u) + \int_{0}^\infty F_1(u + x) dF_1(x) \nonumber \\
& \approx& e^{-1/2} F_1(u) + \sum_{j=0}^M F_1(u + x_j)(F_1(x_j) - F_1(x_{j-1})). \label{e:convPP}
\end{eqnarray}
The first term appears here since $F_2$ has a discontinuity at $x=0$ with
size
\begin{equation*}
F_2(0^+) - F_2(0^-) = e^{-\int_0^{1/2} n(ds)} = e^{-1/2}.
\end{equation*}
This gives an approximation to the CDF of the sum in (\ref{e:breakitupP}).

The CDF and density of $\int_0^1 \sin(2 \pi s) N(ds)$ is  plotted in Figure \ref{f:intsinsNds}.   Similar to Example 2, we see at $x=\pm 1$ two sharp corners in the CDF.  To understand this,
notice that on the interval $(0,1/2)$ the inverses of the function $\krnl(s) = \sin(2 \pi s)$ are given by
\begin{eqnarray*}
\krnl^{-1}(x) &=& \frac{\sin^{-1}(x)}{2 \pi}, \quad \mbox{ for $\krnl$ restricted to $(0,1/4)$}, \\
\krnl^{-1}(x) &=& \frac{\sin^{-1}(1/2 - x)}{2 \pi}, \quad \mbox{ for $\krnl$ restricted to $(1/4,1/2)$}.
\end{eqnarray*}
In both of these cases, $\krnl^{-1} \in C^{0,1/2}(0,1)$ since the
derivative at $x=1$ doesn't exist and
\begin{equation}
\sin^{-1}(1) - \sin^{-1}(x) \sim \sqrt{1 - x}, \ \ \mbox{as} \ x
\rightarrow 1^{-}.
\end{equation}
Thus, Theorem \ref{t:holder} and  Corollary \ref{c:worstcase} imply that $F_1
\in C^{0,1/2}(0,\infty)$, and similar to Example 2, this is evident by
a kink in in $F_1$ at $x=1$.  Since $F_2(x) = 1 - F_1(-x)$, $F_2 \in
C^{0,1/2}(-\infty,0)$, and we expect a kink at $x = -1$.    Because $F_1$ and $F_2$ both have a
discontinuity at $x=0$, the convolution in (\ref{e:convPP}) can be
written in two ways:
\begin{subequations}
\begin{eqnarray}
F(x) &=& \int_{-\infty}^\infty \cdf_1(x - u) d\cdf_2(u)  =  e^{-1/2} F_1(x) + \int_{-\infty}^0 F_1(x - u) dF_2(u)  \label{e:split1}\\
&=& \int_{-\infty}^\infty \cdf_2(x - u) d\cdf_1(u)  =  e^{-1/2} F_2(x) + \int_{0}^\infty F_2(x - u) dF_1(u). \label{e:split2}
    \end{eqnarray}
\end{subequations}
Since $\cdf_1$ has a kink at $x = 1$ and $F_2$ has a kink at $x = -1$,
(\ref{e:split1}) and (\ref{e:split2}) explain why you see
kinks at $x =\pm1 $ in the convolution.

\section{Approximation Lemmas}\label{s:approx}

This section contains the bounds on the approximation errors
$S_{j,\ell}(h)$ and $R_{j,\ell}(\delta)$ defined in Section
\ref{s:error}.   Recall that $S_{j,\ell}(h)$ involves an approximation
of the time derivative $\frac{\partial \cdf}{\partial t^+}$ and $R_{j,\ell}(\delta)$ involves an
approximation of $\cdf(x - g(t),t)$.
The first lemma below bounds $S_{j,\ell}(h)$ in terms of
$h,h^{1+\gamma}$ and $h^2$.   The second lemma bounds
$R_{j,\ell}(\delta)$ in terms of $\delta^{1+\gamma}$.  Here $\gamma$
is the H\"{o}lder exponent of $\krnl^{-1}$.

\begin{lem}\label{l:Slem}
Let $S_{j,\ell}(h)$ be as in (\ref{e:defofS}), and assume $g$ and $n$
satisfy the assumptions in (\ref{e:gnassump}).   Then
\begin{equation}\label{e:Slem}
|S_{j,\ell}(h) | \leq \begin{cases} \displaystyle n_T^{\star} (2 h + n^{\star}_T h^2
  ) \qquad x_j - \krnl(t_{\ell}) \geq 0 > x_j - \krnl(t_{\ell + 1}) \\ 
 \displaystyle  (\frac{L_n}{2}+2) n^{\star \ 2}_T h^2 + \frac{ n^{\star}_T C_{n,T} L_\krnl^\gamma [\krnl^{-1}]_{C^{0,\gamma}(\krnl(0),\krnl(T))} }{1 + \gamma} h^{1 + \gamma}, \quad  x_j - \krnl(t_{\ell + 1}) \geq 0 \\
(\frac{n_T^{\star} L_n}{2} + 2 n^{\star}_T ) h^2 \qquad x_j - \krnl(t_\ell) < 0 \end{cases}
\end{equation}
 where $C_{n,T}$ is defined in Theorem \ref{t:holder}.
\end{lem}

\begin{pf}
Let $ h >0$.  As in the proof of Theorem \ref{t:mainthm}, we have 
\begin{eqnarray}
\cdf(x_j,t_{\ell}+h) - \cdf(x_j,t_{\ell}) &=& \mathbb{E} [
\bone_{(-\infty,x]}(X(t_{\ell}+h)) - \bone_{(-\infty,x]}(X(t_{\ell}))
]  \nonumber \\
&=&\sum_{n=0}^\infty \mathbb{E} [\bone_{(-\infty,x]}(X(t_\ell+h)) -
\bone_{(-\infty,x]}(X(t_\ell)) | N = n] P_n \nonumber \\
&=& \left(  \int_{t_\ell}^{t_\ell+h} ( \cdf(x_j - \krnl(s),t_\ell) - \cdf(x_j,t_\ell)) n(s) ds  \right) e^{-\int_{t_\ell}^{t_\ell+h} n(s)ds} \nonumber   \\
&\ &  + \sum_{n=2}^\infty \mathbb{E} [\bone_{(-\infty,x]}(X(t_\ell+h)) - \bone_{(-\infty,x]}(X(t_\ell)) | N = n] P_n  \label{e:ALemstep1}
\end{eqnarray}
Now, using (\ref{e:mainde}), (\ref{e:ALemstep1}), and the fact that $ 0  \leq \cdf(x_j,t_\ell) \leq 1$, we have
\begin{eqnarray}
|S_{j,\ell}(h) | &=&  | \cdf(x_j,t_\ell+h) - \cdf(x_j,t_\ell) - h
\frac{\partial \cdf}{\partial t^+} | \nonumber \\
&=& | \cdf(x_j,t_\ell+h) - \cdf(x_j,t_\ell) - h(\cdf(x_j - \krnl(t_\ell),t_\ell) - \cdf(x_j,t_\ell))n(t_\ell) | \nonumber \\ 
&\leq& \left| \left(  \int_{t_\ell}^{t_\ell+h} ( \cdf(x_j - \krnl(s),t_\ell) - \cdf(x_j,t_\ell)) n(s) ds  \right)  - h  ( \cdf(x_j - \krnl(t_\ell),t_\ell) - \cdf(x_j,t_\ell)) n(t_\ell)  \right| \nonumber \\
&\ & + \left| \left(  \int_{t_\ell}^{t_\ell+h} ( \cdf(x_j - \krnl(s),t_\ell) - \cdf(x_j,t_\ell)) n(s) ds  \right) (e^{-\int_{t_\ell}^{t_\ell+h} n(s)ds} - 1) \right|  + \sum_{n=2}^\infty P_n  \label{e:ALemstep2} 
\end{eqnarray}
To obtain (\ref{e:Slem}), we consider the three cases separately.
We will use the following elementary inequalities involving $e^{-x}$
for $x \in \mathbb{R}$: 
\begin{subequations}
\begin{eqnarray}
1 - e^{-x}  &\leq&  x \label{e:ebound1}\\ 
1 - e^{-x} - x e^{-x} & \leq&   x^2. \label{e:ebound2}
\end{eqnarray}
\end{subequations}

\medskip

\textbf{Case 1.} $x_j - \krnl(t_{\ell}) \geq 0 > x_j - \krnl(t_{\ell + 1}) $

\medskip

This case contains the ``troublesome'' point, since $ 0 \notin
C(\cdf)$ and Theorem \ref{t:mainthm} does not necessarily apply.
However, in view of (\ref{e:leftversion}) in Remark 2 after the proof of Theorem \ref{t:mainthm}, we have in any case
\begin{equation}
\left|h \frac{\partial \cdf}{\partial t_\ell^{+}} \right| \leq h n^{\star}_T.
\end{equation} 
And, from (\ref{e:ALemstep1}),
\begin{eqnarray*}
|\cdf(x_j,t_\ell+h) - \cdf(x_j,t_\ell)| &\leq& n^{\star}_T h +
\sum_{n=2}^\infty P_n.
\end{eqnarray*}
Using (\ref{e:ebound2}), we have
\begin{eqnarray}
\sum_{n=2}^\infty P_n &=& 1 - P_0 - P_1  = 1 -  e^{-\int_t^{t+h} n(s)
  ds} - \left( \int_t^{t+h} n(s) ds \right)e^{-\int_t^{t+h} n(s) ds}
\nonumber 
\\
 &\leq& \left( \int_t^{t+h} n(s) ds \right)^2 \leq  (h n^{\star}_T)^2. \label{e:lemp3}
\end{eqnarray} 
Thus,  for this case, 
\begin{equation}
|S_{j,\ell}(h) | \leq |\cdf(x_j,t_\ell+h) - \cdf(x_j,t_\ell)| +  \left|h \frac{\partial \cdf}{\partial t^{+}} \right| \leq 2 n^{\star}_T h + n^{\star \ 2}_T h^2,
\end{equation}
which is consistent with (\ref{e:Slem}).

\medskip

\noindent \textbf{Case 2}: $x_j - \krnl(t_{\ell + 1}) \geq 0$.

\medskip

Since $\krnl$ is monotone increasing, $x_j - \krnl(t_\ell) > x_j -
\krnl(t_{\ell + 1}) \geq 0$.  From Theorem \ref{t:holder} and the
Lipschitz assumptions on $n$ and $\krnl$ imply the first term in
absolute values in (\ref{e:ALemstep2}) can be
rewritten and bounded as
\begin{equation*}
\left| \int_{t_\ell}^{t_\ell+h} (\cdf(x_j - \krnl(s),t_\ell) - \cdf(x_j - \krnl(t_\ell),t_\ell))n(s) ds  -  \int_{t_\ell}^{t_\ell+h} (n(t_\ell) - n(s))(\cdf(x_j - \krnl(t_\ell),t_\ell) - \cdf(x_j,t_\ell)) ds  \right|
\end{equation*}
\begin{eqnarray}
&\leq& n^{\star}_T \int_{t_\ell}^{t_\ell+h} |\cdf(x_j -
\krnl(s),t_\ell) - \cdf(x_j - \krnl(t_\ell),t_\ell))| ds +
\int_{t_\ell}^{t_\ell+h} |n(s) - n(t_\ell) | ds. \nonumber \\
&\leq&  n^{\star}_T  C_{n,T} [\krnl^{-1}]_{C^{0,\gamma}(\krnl(0),\krnl(T))}   \int_{t_\ell}^{t_\ell+h} |\krnl(s) -
\krnl(t)|^\gamma ds +
\int_{t_\ell}^{t_\ell+h} |n(s) - n(t_\ell) | ds. \nonumber \\
&\leq& n^{\star}_T C_{n,T} [\krnl^{-1}]_{C^{0,\gamma}(\krnl(0),\krnl(T))} L_{\krnl}^\gamma \int_{t_\ell}^{t_\ell+h} |s - t_\ell|^\gamma ds + L_n \int_{t_\ell}^{t_\ell+h} |s-t_\ell| ds \nonumber \\
&=& n^{\star}_T C_{n,T} [\krnl^{-1}]_{C^{0,\gamma}(\krnl(0),\krnl(T))} L_{\krnl}^\gamma \frac{h^{1+\gamma}}{1 + \gamma} + L_n \frac{h^2}{2} \label{e:lemp1}
\end{eqnarray} 
From (\ref{e:ebound1}) and since $0 \leq |\cdf| \leq 1$, the second
term in (\ref{e:ALemstep2}) is bounded by 
\begin{eqnarray}
\left| \left(  \int_{t_\ell}^{t_\ell+h} ( \cdf(x_j - \krnl(s),t_\ell) - \cdf(x_j,t_\ell)) n(s) ds  \right) (e^{-\int_{t_\ell}^{t_\ell+h} n(s)ds} - 1) \right| &\leq& h n^{\star}_T  \int_{t_\ell}^{t_\ell+h} n(s) ds \nonumber  \\
& \leq & (h n^{\star}_T)^2, \label{e:lemp2}
\end{eqnarray}
Finally, the third term in (\ref{e:ALemstep2}) is bounded by
(\ref{e:lemp3}).   Putting (\ref{e:lemp1}), (\ref{e:lemp2}) and (\ref{e:lemp3}) together, gives
\begin{equation}
|S_{j,\ell}(h)| \leq (\frac{L_n}{2} + 2) n^{\star \ 2}_T h^2 + \frac{ n^{\star}_T C_{n,T}  [\krnl^{-1}]_{C^{0,\gamma}(\krnl(0),\krnl(t))} L_{\krnl}^\gamma }{1 + \gamma} h^{1 + \gamma }
\end{equation}
which agrees with (\ref{e:Slem}) for this case.

\medskip

\textbf{Case 3.}  $x_j - \krnl(t_\ell) < 0$

In this case,  many of the terms in (\ref{e:ALemstep2}) drop out since $\cdf(x,t) = 0$ for $x < 0$ and $\krnl$ is assumed to be non-negative.  Breaking up (\ref{e:ALemstep2}) in a similar way to case 2, we obtain 
\begin{eqnarray*}
|S_{j,\ell}(h) | &\leq&  \int_{t_\ell}^{t_\ell+h} |n(s) - n(t_\ell)| \cdf(x_j,t_\ell) ds + 2 n^{\star \ 2}_T h^2 \\
&\leq & (\frac{n^{\star}_T L_n}{2} + 2 n^{\star}_T )h^2
\end{eqnarray*}
This confirms (\ref{e:Slem}) and finishes the proof.

\end{pf}

\begin{lem}\label{l:Rlem}
Let $R_{j,\ell}(\delta)$ be as in (\ref{e:defofR}), and assume $\krnl$
satisfies the assumptions in (\ref{e:gnassump}).   Then
\begin{equation}
|R_{j,\ell}(\delta)| \leq C_{n,T} [\krnl^{-1}]_{C^{0,\gamma}(\krnl(0),\krnl(T))} \delta^{\gamma},
\end{equation}
where $C_{n,T}$ is defined in Theorem \ref{t:holder}.
\end{lem}

\begin{pf}

By the definitions of $\lambda_{\ell}$ and $k_\ell$, and since $\cdf$
and $\krnl$ are non-decreasing, we have
\begin{align*}
\cdf(x_{j -k_\ell},t_\ell) \  &\leq \  \cdf(x_j - \krnl(t_\ell),t_\ell) \  \leq \  \cdf(x_{j -k_\ell+1},t_\ell) \\
\cdf(x_{j -k_\ell},t_\ell) \  &\leq \   (1- \lambda_\ell) \cdf(x_{j - k_\ell},t_\ell) + \lambda_\ell   \cdf(x_{j - k_\ell+1},t_\ell) \ \leq \  \cdf(x_{j -k_\ell+1},t_\ell)  .
\end{align*}
These inequalities and Theorem \ref{t:holder} now imply
\begin{eqnarray}
|R_{j,\ell}(\delta)| &\leq& |  \cdf(x_{j -k_\ell+1},t_\ell)   - \cdf(x_{j -k_\ell},t_\ell) | \nonumber \\
&\leq& C_{n,T} [\krnl^{-1}]_{C^{0,\gamma}(\krnl(0),\krnl(T))} |x_{j - k_\ell + 1} - x_{j - k_\ell}|^\gamma \nonumber \\
&=& C_{n,T} [\krnl^{-1}]_{C^{0,\gamma}(\krnl(0),\krnl(T))} \delta^\gamma.
\end{eqnarray}
This finishes the proof.

\end{pf}

\section{Guide to Software}\label{s:software}

Software written in MATLAB for computing the CDF of a Poisson integral of the form
(\ref{e:firstpoisint}) is freely available from the authors.
This section will go over installation and use.  

Once the file \texttt{Poisson\_Integral\_GUI.zip} is downloaded and
decompressed, move the folder ``\texttt{Poisson\_Integral\_GUI}'' to
your desired directory.  Launch MATLAB, and enter
\texttt{path(path,'your\_path/Poisson\_Integral\_GUI')},  where \texttt{your\_path} denotes the path to the folder
\texttt{Poisson\_Integral\_GUI}.  You should now be able to use the code
included with the package.

\begin{figure}[h]
\centering
\includegraphics[scale=.4]{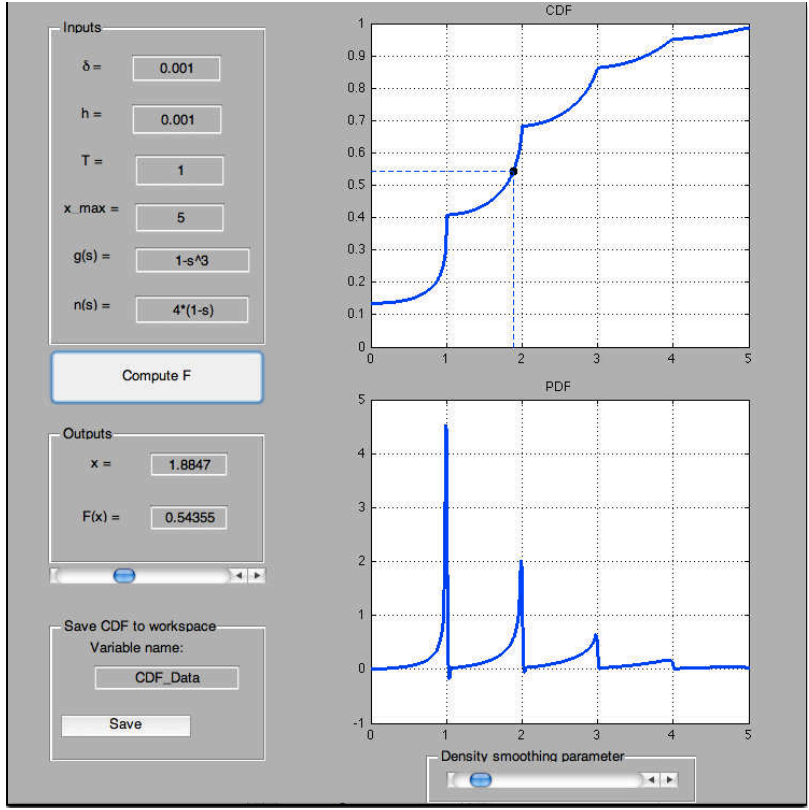}
\caption{ A snapshot of the GUI. } \label{f:GUI}
\end{figure}

To begin, type \texttt{Poisson\_Integral\_GUI} and press return.  This
will open a window similar to that shown in figure \ref{f:GUI}.  From
here you may enter the parameters into the
``Inputs'' panel on the top left.  These are:
\begin{align*}
&\delta \quad &\mbox{Step size in spatial mesh} \\
&h \quad &\mbox{Step size in temporal mesh} \\
&T \quad &\mbox{Upper limit of integration in (\ref{e:firstpoisint})} \\
&x_{\max} \quad &\mbox{Largest value for which CDF is computed} \\
&g(s) \quad &\mbox{Kernel function} \\
&n(s) \quad &\mbox{Density of the control measure}
\end{align*}  
The kernel function and control measure should be entered as
MATLAB expressions in terms of $s$.  For example,  the kernel $g(s) =\sin^2(2\pi
s)$ and Lebesgue control measure $n(s) = 1$ are entered as
\texttt{sin(2 * pi * s)$\wedge$2} and \texttt{1}, respectively.  Only
non-negative kernels and control measures will be accepted.

Once all the parameters are entered, clicking on \texttt{Compute F} will
generate the CDF of the integral $\int_0^T g(s) N(ds)$ in the upper
right plot, as well as a ``quick and dirty'' approximation of the
continuous part of the density in the lower right plot\footnote{This
  requires the Spline Toolbox.}.  The density
is approximated by fitting a smoothing spline to the CDF, then
computing the derivative of this spline.  The slider
under the density plot allows you to adjust the smoothing parameter
used in computation of the density.  

Once a CDF is computed, you may compute specific values of the CDF
using the ``Outputs'' panel on the left.  Simply enter in a desired value
in the $x$ slot, and press return.  Alternatively, you may select an
$x$-value using the slide bar at the bottom of the panel.  You can also enter a
value in the $F(x)$ slot to compute the inverse of the CDF.  

The lower left panel allows you to save the CDF data to your MATLAB workspace.  To do this, enter a
variable name, and click \texttt{save}.  This will create a variable
on your workspace with your chosen name.  This variable is an $M$ by 2
array whose first column contains the $x$ values
$(0,\delta,2\delta,\dots,M\delta)'$ and second column contains the
corresponding CDF values, i.e.~$(F(0),F(\delta),\dots,F(M \delta))'$
(here, $M$ denotes the size of the spatial mesh).   

Finally, if the user wants to avoid using the interface, the code
\texttt{poissonCDF} will work in the command line.  This function
takes in the six parameters listed above and outputs a vector
containing the values of the CDF on the mesh points generated.  The
functions $g$ and $n$ must be entered as function handles.  For
example, if $\delta = 0.01, h= 0.01, T=3, x_{\max} = 4, g(s) = s^2$
and $n(s) = 1$, the call to this function would look like

\medskip

\noindent \texttt{F = poissonCDF(.01,.01,4,3,@(s) s $\wedge$ 2,@(s) 1 );}

\medskip

\noindent The ``\texttt{@}'' symbol is necessary because the input
must be a so-called ``function
handle''.  The statement above will return the vector \texttt{F} which
contains the values of the CDF of the integral $\int_0^3 s^2 N(ds)$
with Lebesgue control measure for the values $x =
.01,.02,\dots,3.99,4$.  For more on this function, read the comments in the code.

\noindent Mark Veillette \& Murad Taqqu \\
\noindent \small Dept. of Mathematics \\
\noindent \small Boston University \\
\noindent \small 111 Cummington St. \\
\noindent \small Boston, MA 02215

\end{document}